\documentclass[preprint,12pt]{elsarticle}
\usepackage{graphicx}

\usepackage[english]{babel}
\usepackage[T1]{fontenc}  
\usepackage[utf8]{inputenc}  
\usepackage{textcomp}  
\usepackage{amsmath, amsfonts, amssymb, bm}
\usepackage[dvipsnames]{xcolor}
\usepackage{comment}
\usepackage[noend]{algpseudocode}
\usepackage{algorithm}
\algnewcommand\algorithmicto{\textbf{to}}
\algnewcommand\algorithmicin{\textbf{in}}
\algnewcommand{\algorithmicand}{\textbf{ and }}
\algnewcommand{\algorithmicor}{\textbf{ or }}
\algnewcommand{\algorithmicnot}{\textbf{not }}
\algnewcommand\algorithmicforeach{\textbf{for each}}
\algrenewtext{For}[3]{\algorithmicfor\ #1 $\gets$ #2\ \algorithmicto\ #3\ \algorithmicdo}
\algdef{S}[FOR]{ForEach}[2]{\algorithmicforeach\ #1\ \algorithmicin\ #2\ \algorithmicdo}
\algnewcommand{\OR}{\algorithmicor}
\algnewcommand{\AND}{\algorithmicand}
\algnewcommand{\NOT}{\algorithmicnot}

\usepackage{hyperref}
\hypersetup{
    colorlinks = true,
    linkcolor = blue,
    linkbordercolor = {white}
}

\renewcommand{\vec}[1]{\boldsymbol{#1}}
\newcommand{\mt}[1]{\bm{#1}}
\newcommand{\R}{\mathbb{R}}
\newcommand{\re}{\mathrm{Re}}

\newcommand{\norm}[1]{\|#1\|}
\newcommand{\abs}[1]{\left\lvert#1\right\rvert}

\graphicspath{}

\begin{document}
\begin{frontmatter}

\title{Adaptive hyperviscosity stabilisation for the RBF-FD method in solving advection-dominated transport equations} 

\author[1,2]{Miha Rot\corref{cor1}}
\ead{miha.rot@ijs.si}
\author[1,3]{\v{Z}iga Vaupoti\v{c}}
\author[1,3]{Andrej Kolar-Po\v{z}un}
\author[1]{Gregor Kosec}
\ead{gregor.kosec@ijs.si}

\affiliation[1]{organization={Parallel and Distributed Systems Laboratory, Jo\v{z}ef Stefan Institute},
                addressline={Jamova cesta 39},
                postcode={1000},
                postcodesep={},
                city={Ljubljana},
                 country={Slovenia}}

\affiliation[2]{organization={Jozef Stefan International Postgraduate School},
                addressline={Jamova cesta 39},
                postcode={1000},
                postcodesep={},
                city={Ljubljana},
                country={Slovenia}}

\affiliation[3]{organization={Faculty of Mathematics and Physics, University of Ljubljana},
                addressline={Jadranska ulica 19},
                postcode={1000},
                postcodesep={},
                city={Ljubljana},
                country={Slovenia}}
 
\cortext[cor1]{Corresponding author}

\begin{abstract}

This paper presents an adaptive hyperviscosity stabilisation procedure for the Radial Basis Function-generated Finite Difference (RBF-FD) method, aimed at solving linear and non-linear advection-dominated transport equations on domains without a boundary. The approach employs an equation independent algorithm that adaptively determines the hyperviscosity constant based on the spectral radius of the RBF-FD evolution matrix. The proposed procedure supports general node layouts and is not tailored for specific equations, avoiding the limitations of empirical tuning and von Neumann-based estimates. To reduce computational cost, it is shown that lower monomial augmentation in the approximation of the hyperviscosity operator can still ensure consistent stabilisation, enabling the use of smaller stencils and improving overall efficiency. A hybrid strategy employing different spline orders for the advection and hyperviscosity operators is also implemented to enhance stability. The method is evaluated on pure linear advection and non-linear Burgers' equation, demonstrating stable performance with limited numerical dissipation. The two main contributions are: (1) a general hyperviscosity RBF-FD solution procedure demonstrated on both linear and non-linear advection-dominated problems, and (2) an in-depth analysis of the behaviour of hyperviscosity within the RBF-FD framework, addressing the interplay between key free parameters and their influence on numerical results.
\end{abstract}


\end{frontmatter}

\section{Introduction}
Transport Partial Differential Equations (PDE) arise from fundamental conservation principles and provide a mathematical framework to describe how conserved quantities move and evolve over time and space. Examples of such PDEs include the advection equation for the transport of a quantity by a moving fluid, the diffusion equation for describing the spreading of a substance due to random motion, the Navier-Stokes equation governing fluid motion through momentum conservation, and its simplified variant, Burgers' equation, capturing key non-linear and viscous effects while omitting pressure, to name a few.
Solving these equations is crucial for understanding physical processes in various fields ranging from pure scientific research, to environmental modelling and real-world engineering applications. Transport equations are typically solved numerically due to their complexity. Nevertheless, even numerical simulations, especially in advection-dominated regimes, can be challenging. If not managed correctly, they may produce non-physical oscillations, which can potentially lead to divergence~\cite{sirca_computational_2018}.

In most numerical simulations, mesh-based methods such as the Finite Volume Method (FVM), or the Finite Element Method (FEM) are used. However, due to meshing limitations an alternative class of meshless methods emerged in 1970s~\cite{Liu_2002}. The fundamental distinction between traditional mesh-based methods and meshless methods lies in how they establish relationships between computational nodes. Unlike mesh-based approaches, meshless methods define these relationships solely based on relative node positions, allowing them to operate on scattered nodes without relying on any predefined structure or mesh~\cite{nguyen_meshless_2008}. Although achieving a stable meshless approximation requires scattered nodes to follow specific guidelines~\cite{Slak2019,van_der_sande_fast_2021}, the process of positioning these nodes is far less complex than traditional meshing~\cite{Zienkiewicz_Taylor_Zhu_2005}. This makes meshless methods highly flexible for handling complex geometries, including those supplied by Computer-Aided Design (CAD) models~\cite{Mirfatah_Boroomand_Soleimanifar_2019,duh_discretization_2024}, as well as for applications in h-adaptivity~\cite{Reeger_2024,Oanh_Davydov_Phu_2017,Slak_Kosec_2019b} and moving boundary problems~\cite{ortega2013meshless,jancic_meshless_2024}. Scattered node positioning algorithms can also be implemented in a dimension-independent manner~\cite{Slak2019,Shankar_Kirby_Fogelson_2018} that can be combined with the elegant formulation of meshless approximation techniques to facilitate solution procedures for high-dimensional problems~\cite{Jancic_Slak_Kosec_2021}. 

Among emerging meshless methods, the Radial Basis Function-generated Finite Differences (RBF-FD)\cite{Tolstykh_Shirobokov_2003} has in the last two decades gained significant attention for its adaptability and robustness in solving various PDEs. RBF-FD most commonly employs polyharmonic spline (PHS) radial basis functions (RBFs) augmented with monomials \cite{Bayona_Flyer_Fornberg_Barnett_2017,bayona2019insight} for differential operator approximation, which helps prevent stagnation errors and allows control over the convergence rate. Besides improved stability, the flexibility in the order of approximation has recently led to the development of the first truly meshless hp-adaptive method~\cite{Jancic_Kosec_2024}, high-order solutions for various problems~\cite{Shahane_Radhakrishnan_Vanka_2021,nielsen2024high,Jancic_Slak_Kosec_2021}, and the implementation of spatially variable operator approximations to improve computational performance~\cite{rot2024spatially}.

Nonetheless, despite the potential advantages of the RBF-FD method, its differentiation matrices inherently contain spurious eigenvalues~\cite{Platte2006,Fornberg2011}. This effect remained theoretically unexplained until recently, when Tominec et al. ~\cite{Tominec2025} provided a mathematical foundation for it within the framework of energy methods, a common approach in finite element analysis. The authors demonstrated that instabilities stem from integration errors and discontinuities in global cardinal functions, which were already discussed in~\cite{Tominec2021}. In that work, the authors introduced an oversampling technique, a method for stabilising RBF partition of unity method (RBF-PUM)~\cite{Larsson2017}
and Kansa formulations~\cite{Tominec2025,Martel2016}, where global cardinal functions do not exhibit jumps. On the other hand, the jumps in global cardinal functions present in the RBF-FD method can be reduced by increasing the stencil size~\cite{Tominec2025}, a conditionally useful approach due to its high computational cost.

In practice, dissipative numerical schemes are used for stabilisation. 
Arguably the most widely used dissipative stabilisation technique is upwind~\cite{golbabai2019analysis}. Upwind is a directional biasing approach in which, when computing numerical derivatives, greater weight is assigned to values from the direction where information is coming from. Upwind schemes are widely used in computational fluid dynamics and numerical simulations where convective transport is dominant~\cite{lyra2000review}.
Nonetheless, while upwind stabilises numerical simulations that would otherwise diverge, it also introduces significant artificial diffusion~\cite{lyra2000review}, which consequently smears out sharp features in the final solution~\cite{kosec_simulation_2014}.

In an alternative stabilisation technique, referred to as hyperviscosity, a higher-order diffusion term, i.e. fourth or higher order derivatives, is added to the governing equations. Such hyperviscosity operator selectively damps small-scale oscillations while preserving large-scale structures, i.e. it controls high-frequency instabilities. The concept of a hyperviscosity operator dates back to the 1990s \cite{Ma1998,Boyd1994} and was later adapted to RBF-FD methods~\cite{Fornberg2011}, where authors demonstrated its effectiveness in shifting the eigenvalues of the differentiation matrix for convective PDEs to a stable region. In the context of RBF-FD, hyperviscosity has been recently also demonstrated on advection-diffusion PDEs~\cite{shankarhyperviscosity} and advection-diffusion-reaction PDEs on manifolds~\cite{shankar2019}.

The main problem with the hyperviscosity scheme is its inherent constant, typically denoted by $\gamma$, which must be tuned specifically for each case. Although an empirical estimate for this constant was proposed~\cite{Flyer2016NS,Flyer2012}, user input was still required to determine its precise value. It was not until recently that the first attempt to theoretically determine an appropriate $\gamma$ was introduced by Shankar et al.~\cite{shankarhyperviscosity}, using a quasi-analytical approach based on von Neumann analysis. However, in addition to standard limitations of von Neumann analysis, such as its validity only for periodic or infinite domains, linear PDEs, uniform spatial discretisation, and its disregard for boundary effects, the proposed approach also requires an estimate of an additional growth term to account for scattered nodes and RBF-FD interpolation. This term is not well-defined, which may lead to inaccurate results and, consequently, an ill-determined $\gamma$. 

This paper advances the applicability of hyperviscosity by refining the understanding and determination of hyperviscosity constants. The primary goal is to develop a PDE-independent approach for optimising the parameters of the hyperviscosity scheme on domains without a boundary and to enhance the understanding of their impact on accuracy and stability. We demonstrate that an appropriate hyperviscosity constant can be estimated through the eigenvalues of the evolution matrix, as indicated by Fornberg et al.~\cite{Fornberg2011}, by using a relatively simple iterative algorithm. The proposed approach is not restricted by the type of PDE, domain geometry, RBF-FD setup, or the layout of scattered nodes.

Another challenge of hyperviscosity is its high computational cost, as increasing its order significantly enlarges the required stencil. For example, a stable computation of fourth derivatives, as required by the second-order hyperviscosity operator in $2$D, necessitates augmentation with fourth-order monomials and a corresponding stencil of $30$ support nodes~\cite{Bayona_Flyer_Fornberg_Barnett_2017}. Here, we demonstrate that the hyperviscosity operator exhibits consistent behaviour even by using a smaller monomial augmentation in constructing the approximation. This finding suggests that computational costs can be reduced by employing only second-order augmentation, regardless of the order of the hyperviscosity operator.
The developed methodology is demonstrated on a linear advection-dominated problem~\cite{Tominec2022,shankarhyperviscosity} and a non-linear Burgers' problem~\cite{Fletcher1983,Zhu2010}, where we clearly illustrate its efficiency. 

The paper is organised in the following order: In \autoref{sec:rbffd} we provide a description of the RBF-FD method and use it to discretise linear and non-linear Cauchy problems in the matrix-vector format. In \autoref{sec:hyperviscosity} we construct a hyperviscosity operator and theoretically analyse its parameters. The hyperviscosity scheme and parameters are first numerically verified on simple linear advection case in \autoref{sec:linearadvection}. These results are extended to include the non-linear Burgers' equation in \autoref{sec:burgers}. We finish the paper with some remarks on the future work in \autoref{sec:discussion}.

\section{The RBF-FD method}
\label{sec:rbffd}
First, we discuss the RBF-FD approach to discretising differential operators on scattered nodes spanning the domain $\Omega$ and therefore constructing the appropriate differentiation matrices.

\subsection{RBF approximation}
\label{sec:rbffd_approx}
The domain $\Omega \subseteq \mathbb{R}^d$ is first discretised using a set of quasi-uniformly distributed scattered nodes $X=\{\vec{x}_i\}_{i=1}^N \subset \Omega$ that are generated using a dedicated advancing front algorithm~\cite{Slak2019} with a minimal spacing quarantee. We use $h = \sup_{x_i \in X} \min_{x_j \in X \setminus \{x_i\}} \| x_i - x_j\|$ to denote the largest distance between first neighbours. Each node $\vec{x}_i \in X$ is assigned a subset $X^{(i)} \subseteq X$ of the closest nodes $n$ (commonly referred to as a stencil or support domain) that are used for construction of a local interpolation in the vicinity of that node. Using RBFs $\phi_i(\vec{x}) = \phi (\norm{\vec{x}_i - \vec{x}}), i=1,\dots, n$ and a multivariate monomial basis $\{m_j\}_{j=1}^q $, the local interpolation reads
\begin{equation}
    \begin{gathered}
    u_h(\vec{x}) = \sum_{i=1}^n w_i \phi_i(\vec{x}) + \sum_{i=1}^q \beta_i m_i(\vec{x}), \\ \text{subject to: } \, \sum_{i=1}^n w_i m_j(\vec{x}_i) = 0 \quad j= 1,\dots, q,
    \end{gathered}
    \label{eq:pointwise_rbf_interp}
\end{equation}
where $w_i$ are the RBF interpolation weights and $\beta_i$ are the Lagrange multipliers. To obtain the weights a linear system is assembled enforcing the pointwise interpolation in equation~\eqref{eq:pointwise_rbf_interp} for each point in the stencil:
\begin{align}
     \begin{bmatrix}
\mt{A} & \mt{P} \\ \mt{P}^T & 0 
\end{bmatrix}    \begin{bmatrix}
\vec{w} \\ \vec{\beta} 
\end{bmatrix} = \begin{bmatrix}
    \vec{u_h} \\ 0
\end{bmatrix},
\label{eq:local_interp_system}
\end{align}
where $\vec{w} = [w_i]_{i=1}^n$ is the interpolation weight vector, $\vec{\beta} = [\beta_i]_{i=1}^q $ is the Lagrange multipliers vector, $\vec{u_h} = [u_h(\vec{x}_i)]_{i=1}^n$ is a vector of stencil nodal values, and $\mt{A}_{i,j} = \phi_i(\vec{x}_j)$ and $\mt{P}_{i,j} = m_j(\vec{x}_i)$ are interpolation matrices. In this paper, we use conditionally positive definite RBFs, more specifically, polyharmonic splines $\phi(r)=r^k$, where $r = \norm{\vec{x}_j-\vec{x}}_2$ is the Euclidean distance to the centre and $k \in \mathbb{N}$ is odd. Monomials are added to the local interpolation \eqref{eq:pointwise_rbf_interp} to ensure the unisolvency of the local interpolation system \eqref{eq:local_interp_system} for conditionally positive definite RBFs \cite{bayona2019,Fasshauer2007} and to improve the convergence characteristics \cite{Jancic_Slak_Kosec_2021,LeBorne2023}. To prove unisolvency, we require a monomial tensor basis $\{m_j\}_{j=1}^q $ of size $q=\binom{m+d}{m}$ to span the polynomial space of order $m$, where $m \geq \lceil \frac{k}{2} \rceil - 1$ \cite{Fasshauer2007}.
The recommended stencil size $n$ for $2$D systems is $n= 2q$~\cite{bayona2019}. 

In \autoref{sec:consistency}, where we discuss the consistency of the hyperviscosity scheme, we will use the concept of global cardinal functions. To approximate the linear differential operator $\mathcal{L}$ acting on function $u: \Omega \to \mathbb{R}^d$, we use the global cardinal functions $\Psi_i: \Omega \to \mathbb{R}, i=1,\dots, N$ with the following ansatz
\begin{align}
	u_h (\vec{x}) = \sum_{i=1}^N \Psi_i(\vec{x}) u_h (\vec{x}_i) \qquad  (\mathcal{L}u_h)(\vec{x}) = \sum_{i=1}^N \mathcal{L} \Psi_i(\vec{x}) u_h (\vec{x}_i),
	\label{eq:lwu}
\end{align}
where $u_h(\vec{x}_1), \dots, u_h(\vec{x}_N) $ are the degrees of freedom. 
The global cardinal functions are constructed by blending the local interpolation systems for every stencil in the domain. A local interpolation system for node $i$ is (re)introduced via local cardinal functions $\psi_j^{(i)} : \Omega \to \mathbb{R}$ and values $u_h(\vec{x}_j)$ in stencil nodes $\vec{x}_j \in X^{(i)}$
\begin{align}
	u_h^{(i)}(\vec{x}) = \sum_{x_j \in X^{(i)}} \psi_j^{(i)}(\vec{x}) u_h(\vec{x}_j) \qquad (\mathcal{L}u_h^{(i)})(\vec{x}) = \sum_{x_j \in X^{(i)}} \mathcal{L} \psi_j^{(i)}(\vec{x}) u_h(\vec{x}_j).
	\label{eq:local_cardinal_functions}
\end{align} 
The equation \eqref{eq:pointwise_rbf_interp} can now be rewritten using the local cardinal functions \eqref{eq:local_cardinal_functions} as

\begin{equation}
\begin{aligned}
    u_h^{(i)}(\vec{x}) &= \begin{bmatrix}
        \vec{\phi}^{(i)}(\vec{x}) & \vec{m}^{(i)}(\vec{x})
    \end{bmatrix} \begin{bmatrix}
        \vec{w}^{(i)} \\ \vec{\beta}^{(i)}
    \end{bmatrix} \\
    &=  \begin{bmatrix}
        \vec{\phi}^{(i)}(\vec{x}) & \vec{m}^{(i)}(\vec{x})
    \end{bmatrix}
    \begin{bmatrix}
        \mt{A}^{(i)} & \mt{P}^{(i)} \\ \left(\mt{P}^{(i)}\right)^T & 0 
    \end{bmatrix}^{-1}
    \begin{bmatrix}
        \vec{u}_h^{(i)} \\ 0
\end{bmatrix} \\&=\begin{bmatrix}
    \vec{\psi}^{(i)}(\vec{x}) & \vec\xi^{(i)}(\vec{x})
\end{bmatrix} \begin{bmatrix}
    \vec{u}_h^{(i)} \\ 0
\end{bmatrix},
\end{aligned}
\end{equation}
where $\vec{\psi}^{(i)}(\vec{x}) = [\psi_j^{(i)} (\vec{x})]_{\vec{x}_j \in X^{(i)}}$ is the local cardinal function vector.

\subsection{Spatial discretisation of PDEs}
To construct the global cardinal functions, for every computational point in the domain we merge all the local cardinal functions used for interpolating that point in local stencils
\begin{align}
    \Psi_i(\vec{x}) = \begin{cases}
        \psi_i^{(\sigma(\vec{x}))}(\vec{x}), \quad &\vec{x}_i \in X^{(\sigma(\vec{x}))}\\
        0, \quad &\vec{x}_i \not \in X^{(\sigma(\vec{x}))}
    \end{cases},
\end{align}
where the function $\sigma(\vec{x})$ returns the index of the closest computational point to any point $\vec{x} \in \Omega$. 
Both global and local cardinal functions have the Kronecker delta property
\begin{align}
    u_h (\vec{x}_j) = \sum_{i=1}^N \Psi_i (\vec{x}_j) u_h(\vec{x}_i) = \sum_{i=1}^N \delta_{i,j} u_h(\vec{x}_i).
\end{align}
 This ensures that the evaluation matrix for sampling the pointwise values of the function is an identity matrix. This also means that the mass matrix does not need to be inverted, which is a considerable overhead in the finite element method. To construct the global differentiation matrix, we assemble the evaluated derivatives of global cardinal functions at every $\vec{x}_i \in X$ in an $N \times N$ matrix. In practice, the evaluated derivatives of cardinal functions are calculated by solving the local interpolation system~\eqref{eq:local_interp_system} with $\begin{bmatrix}\mathcal{L}\vec{\phi}(\vec{x}_i) & \mathcal{L}\vec{m}(\vec{x}_i)\end{bmatrix}$ as the right hand side~\cite{LeBorne2023}.
 
 To put the above-mentioned into perspective, we consider the semidiscrete scheme. Here, we are especially interested in solving linear and non-linear advection, however, for the sake of generality, we consider a Cauchy problem on a domain $\Omega$ without a boundary. Specifically,
 find $u: \Omega \times [0, T] \to \mathbb{R}^{m}$ such that
\begin{align}
\begin{cases}
     \partial_t u(\vec{x},t) + \mathcal{N}(u(\vec{x},t)) = 0 &\quad \text{in $\Omega \times (0, T]$}\\
     u(\vec{x}, 0) = u_0(\vec x) &\quad  \text{in $\Omega$},
     \label{eq:conservation_law}
\end{cases}
 \end{align}
where $\mathcal{N}$ is a (non)linear operator. If $\mathcal{N}$ is nonlinear, we work with its linearisation $\mathcal{L}$. To spatially discretise the equation \eqref{eq:conservation_law} we employ the ansatz \eqref{eq:lwu}  for each point $\vec{x} \in X$ and obtain a set of $N$ equations
\begin{equation}
\begin{aligned}
    \partial_t \sum_{i=1}^N \Psi_i (\vec{x}_1) u_h(\vec{x}_i, t) &= -\sum_{i=1}^N \mathcal{L} \Psi_i (\vec{x}_1) u_h(\vec{x}_i, t), \\
    & \vdots \\
    \partial_t \sum_{i=1}^N \Psi_i (\vec{x}_N) u_h(\vec{x}_i, t) &= -\sum_{i=1}^N  \mathcal{L} \Psi_i (\vec{x}_N) u_h(\vec{x}_i, t). \,
\end{aligned}
\end{equation}
We introduce the differentiation matrix
\begin{equation}
    (\mt{D}_h^\mathcal{L})_{j,i} = \mathcal{L} \Psi_i (x_j)
    \label{eq:differentiation_matrix}
\end{equation}
to rewrite the obtained set of equations in the vector-matrix form
\begin{equation}
    \partial_t \vec{u_h}(t) = -\mt{D}_h \vec{u_h}(t),
    \label{eq:semi_discrete_scheme}
\end{equation}
where $\vec{u_h} (t) = [ u_h(\bm x_i, t)]_{i=1}^N$.To integrate the obtained system of differential equations \eqref{eq:semi_discrete_scheme}, we primarily use the backwards Euler scheme. This particular scheme will be useful due to its fundamental unconditional stability. Given a timestep size of $\Delta t > 0$, we discretise the time domain and obtain
\begin{align}
    \vec{u_h}(t_{n+1}) = \vec{u_h}(t_n) - \Delta t \mt{D}_h \vec{u_h}(t_{n+1}).
    \label{eq:implicit_euler}
\end{align}
In practice, when solving the system \eqref{eq:implicit_euler}, we rewrite it as $(I+\Delta t \mt{D}_h) \vec{u_h}(t_{n+1}) = \vec{u_h} (t_{n})$ and solve the system via LU-decomposition. In the later section, we will also enforce Dirichlet boundary conditions. This is done by replacing the offending rows and columns of the boundary points with the identity matrix.

\section{Hyperviscosity stabilisation}
\label{sec:hyperviscosity}
In this section, we discuss the hyperviscosity stabilisation scheme.
As previously stated RBF-FD differentiation matrices tend to have spurious eigenvalues. Put simply, the RBF-FD solution exhibits numerical artefacts that may lead to instability and divergence of the time-dependent systems. To address this issue, an artificial higher-order Laplacian operator is added to the equation \eqref{eq:conservation_law}~\cite{shankar2019,Ma1998,shankarhyperviscosity,Vaupotic2024}
\begin{equation}
      \partial_t u(\vec x,t) + \mathcal{N} (u(\vec{x}, t)) = (-1)^{\alpha + 1} \gamma \Delta^{\alpha} u(\vec x,t),
      \label{eq:hipscheme}
\end{equation}
or equivalently in the semi-discrete scheme \eqref{eq:semi_discrete_scheme},
\begin{equation}
     \partial_t \vec{u_h}(t) = -\mt{D}_h  \vec{u_h}(t) + (-1)^{\alpha+ 1} \gamma \mt{D}_h^{hip} \vec{u_h}(t).
     \label{eq:semidiscrete_hyp}
\end{equation}
This operator is of particular interest to meshless methods, as it allows us to control the stability of the scheme without relying on a structured mesh. Additionally, the higher order Laplace operator allows us the preserve the order with which our numerical approximation converges to the right solution.

However, as previously discussed, there is no clear consensus in the literature regarding the selection of the constant $\gamma$ and the order $\alpha$. In the following sections, we will discuss the importance of those two parameters and how to select them to satisfy the conditions of the Lax-equivalence theorem.
The stability of the scheme is achieved by correctly choosing $\gamma$ and $\alpha$ terms, which is described in \autoref{sec:stability_gamma}. On the other hand, the scheme is consistent if the RBF-FD parameters are correctly selected; this, along with the optimality of the parameters, is described in \autoref{sec:stability_gamma} and \autoref{sec:parametrs}.

Throughout this section, we will support the discussion with numerical experiments on the stabilised advection operator
\begin{equation}
	\hat{\mt{D}}_h=-\mt{D}_h+ (-1)^{\alpha +1} \gamma \mt{D}_h^{hip}
	\label{eq:advection_op}
\end{equation}
or in some cases on the pure advection operator $\mt{D}_h$.

\subsection{Note on the parameter $\gamma$ and how to iteratively determine it}
\label{sec:stability_gamma}

The estimation of the gamma term has been the subject of numerous papers \cite{Flyer2012,Flyer2016NS,shankar2019,shankarhyperviscosity,Tominec2022} that provided various estimates. Authors of~\cite{Flyer2012} introduced additional scaling that is also used throughout this paper
\begin{equation}
	\gamma(c) = c h^{2\alpha}, 
	\label{eq:c}
\end{equation}
and provided the estimate $c \in [0.1, 1]$ for the Navier-Stokes equation; it is \textbf{assumed} that $c$ is independent or very loosely dependent on $h$. This estimate was later discussed in \cite{Vaupotic2024}. In \cite{Flyer2016} another estimate $\gamma = 2^{-6}  h^{6}$ was given and later generalised by Shankar as  $\gamma = 2^{-2\alpha} h^{2\alpha}$ \cite{shankarhyperviscosity}. Both estimates relied on the Nyquist cut-off wavenumber $k \approx 2h^{-1}$ or what is most commonly referred to as the grid resolution. This general estimate is sometimes sufficient as we aim to target high-order modes with wavenumbers close to the grid resolution. Unfortunately, it only works for simple PDEs and can only be utilised under the assumption that the node set reproduces a grid-like structure. In our case, we primarily choose $h^{2\alpha}$ scaling of the $\gamma$ term to be consistent with the Fourier theory for the hyperviscosity operator \cite{Ma1998,Ma2006}.

\subsubsection{Stability criterion}
To choose $c$, we will introduce a more general approach that relies purely on the eigenvalue decomposition of the discrete evolution matrix. For studying stability it is natural to consider the following weighted $L^2$-like norm~\cite{Tominec2025}
\begin{align}
    \| u \|^2_h =  \frac{|\Omega|}{N}\sum_{i=1}^N  u(x_i) u(x_i) \quad \| \mt{A} \|_h  = \sup_{\norm{x} \ne 0}\frac{\norm{\mt{A}x}_h}{\| x \|_h} = \| \mt{A} \|_2,
\end{align}
where $ {| \Omega| }/{N} \approx h^2$. We must show that our solution is stable with respect to this norm\footnote{Stability in this norm guarantees weak convergence in $L^2$.}.
Firstly, we construct a discrete RBF-FD evolution matrix $\mt{G}_h \in \R^{N\times N}$, that is
\begin{align}
    \vec{u_h} (t_{n+1}) = \mt{G}_h \vec{u_h} ( t_n).
    \label{eq:evolution_matrix}
\end{align}
Note that the hyperviscosity term is directly included in the evolution matrix. For the implicit Euler method, this matrix is
\begin{align}
    \vec{u_h} (t_{n+1}) = \mt{G}_h \vec{u_h}(t_n) = \left[ I - \Delta t \hat{\mt{D}}_h \right]^{-1} \vec{u_h} ( t_{n}),
    \label{eq:dynamical_system_hv}
\end{align}
as discussed regarding the system~\eqref{eq:implicit_euler}. 

The evolution matrix can be constructed for any one-step scheme. By linearity of $\mathcal{L}$, we express the matrix as $\mt {G}_h = \sum_{m=0}^s \xi_m \left(\Delta t \hat{\mt{D}}_h\right)^m$ with scheme specific coefficients $( \xi_i)_{i=0}^s$\footnote{As an example, the coefficients for an explicit linear $s$-stage Runge-Kutta scheme take the form $\xi_m = \frac{1}{m!}$.}.  An evolution matrix can similarly be defined for a linear multi-step method. This is achieved by introducing a lifted state $\bm U_h^{n} = ( \bm {u_h} (t_n), \dots,  \bm {u_h} (t_{n-k+1}))$ and assembling a block matrix from the method's update formulas. While the derivation that follows remains theoretically sound for this type of method, it becomes practically nonviable due to the significant size of the resulting matrix.

Recall that a discretisation of a linear Cauchy problem \eqref{eq:dynamical_system_hv} is said to be \emph{(Lax)}-stable if there exists a constant $C >0$ independent of $h$, such that for all $n \in \mathbb{N}$ the following holds
    \begin{align}
        \norm{\vec{u_h}(t_n)}_h \le C \norm{\vec{u_h} (0)}_h.
    \end{align}
This notion of stability is weaker than the more appropriate energy or entropy stability. However, it still ensures continuous dependence on the initial data. First notice that $\vec{u_h}(t_n) = \mt{G}^n_h \vec{u_h}(0)$.
Assume now that $\mt{G}_h$ is diagonalizable. Then there exists a matrix $\mt{U}_h \in \mathbb{C}^{N \times N}$ and a diagonal matrix $\mt{D} \in \mathbb{C}^{N \times N}$ such that $\mt{G}_h = \mt{U}_h \mt{D}_h \mt{U}_h^{-1}$, this allows us to estimate
\begin{align}
    \| \mt{G}_h^n \vec{u_h}(0) \|_h = \| \mt{U}_h \mt{D}_h^n \mt{U}_h^{-1} \vec{u_h}(0) \|_h \le \norm{\mt{U}_h}_h \norm{\mt{D}_h^n}_h \norm{\mt{U}_h^{-1}}_h \norm{\vec{u_h}(0)}_h.
\end{align}
We denote the condition number by $\kappa (\mt{U}_h) = \norm{\mt{U}_h}_h  \norm{\mt{U}_h^{-1}}_h$ and let $\rho: \mathbb{C}^{N\times N} \to \mathbb{R}^+$ be the spectral radius, i.e. $\rho$ returns the maximal magnitude of the eigenvalues of the matrix. It follows that
\begin{align}
     \| \mt{G}_h^n \vec{u_h}(0) \|_h \le \kappa (\mt{U}_h) \max_{i} \abs{\lambda^n_i} \norm{\vec{u_h}(0)}_h =  \kappa (\mt{U}_h) \rho(\mt{G}_h)^n \norm{\vec{u_h}(0)}_h.
\end{align}
If $\rho(\mt{G}_h) \le 1$, then the system satisfies the stability condition for a fixed $h$. To extend this bound to the Lax-stability criteria, the condition number $\kappa(\mt{U}_h)$ must be uniformly bounded for a sufficiently small $h$. Furthermore, if $\kappa(\mt{U}_h)$ is not close to $1$ then nonphysical growth (and decay) is still possible; hence, we must further assume that $\kappa(\mt{U}_h)$ is close to $1$. We note that RBF-FD matrices are typically not normal and therefore we have no control over the condition number of this matrix.
A similar condition in terms of the spectral radius can be made for the case where $\mt{G}_h$ is not diagonalizable. In this case, we can work with its Jordan canonical form $\mt{G}_h = \mt{U}_h^\prime \mt{J}_h \mt{U}_h^{\prime -1}$ and again estimate
\begin{align}
    \| \mt{G}_h^n \vec{u_h}(0) \|_h \le \kappa (\mt{U}^\prime_h) \norm{\mt{J}^n_h}_h \norm{\vec{u_h}(0)}_h \le \kappa (\mt{U}^\prime_h) \norm{\mt{J}^n_h}_F \norm{\vec{u_h}(0)}_h.
\end{align}
We have used $\norm{\cdot}_2 \le \norm{\cdot}_F$, where $\norm{\cdot}_F$ is the Frobenius norm, which is bounded for all $n$ if $\rho(\mt{G}_h) < 1$ (technically, eigenvalues that correspond to $1 \times 1$ Jordan blocks may have unit modulus)~\cite[Theorem 3.2.5.2.]{Horn2012}. While this guarantees \emph{Lax}-stability, the upper bound does not admit a nice expression in terms of $\rho(\mt{G}_h)$, unlike in the diagonalizable case. In particular, this implies that we do not have any control over the potential nonphysical growth of the solution. Note that from a numerical point of view, this may not be an issue as a generic matrix tends to be diagonalizable. Furthermore, due to round-off errors the conditions $\rho(\mt{G}_h) < 1$ and $\rho(\mt{G}_h) \le 1$ are practically equivalent, justifying the usage of $\rho(\mt{G}_h) \le 1$ as a stability condition for a general case.

\subsubsection{Stabilisation}
\begin{figure}[h!]
	\centering
	\includegraphics[width=1\linewidth]{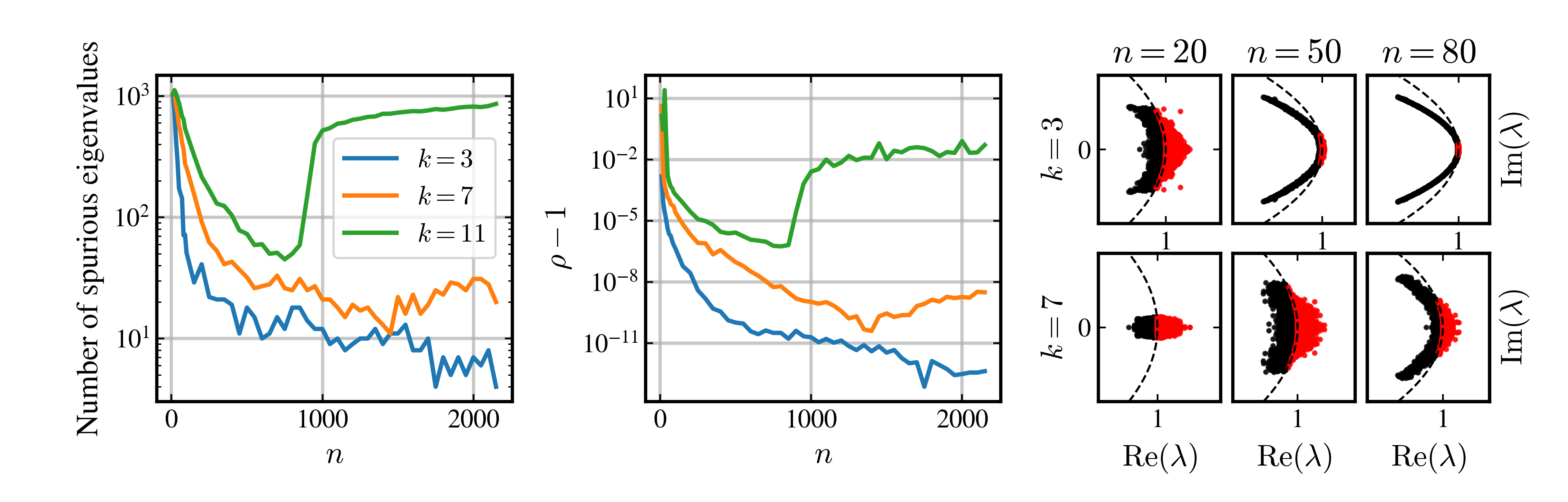}
	\caption{
        Analysis of the impact of stencil size $n$ on the stability of the implicit Euler evolution matrix $\mt{G}_h$~\eqref{eq:evolution_matrix} for an unstabilised linear advection case with $h=0.02$, $\Delta t = 10^{-4}$, $m=2$ and different PHS orders $k$. All eigenvalues of the matrix are computed with a dense decomposition.
        The left graph displays the number of spurious, i.e. with a magnitude above 1, eigenvalues as a function of $n$. The middle graph displays the spectral radius $\rho(\mt{G}_h)$ as a function of $n$. The right cluster of subgraphs displays the actual eigenvalues on the complex plane for a selection of $k$ and $n$. Dashed line marks the unit circle and the spurious eigenvalues outside are coloured red.
}
	\label{fig:stencil_size}
\end{figure}

Now that we have established a condition for stability we can examine how it relates to stabilisation efforts. Before introducing hyperviscosity, we examine a simpler stabilisation approach. Increasing stencil size $n$ used for the approximation is known~\cite{Tominec2025, Shankar2018L} to improve stability and this is reflected in the results displayed in~\autoref{fig:stencil_size}. We can see from the left graph that while the number of spurious eigenvalues with magnitude above 1 decreases drastically as $n$ increases, the spectral radius shown in the middle graph never actually falls into the stable $\rho(\mt{G}_h) \le 1$ regime. Although increasing $n$ improves stability, it also drastically increases the computational costs.
Additionally, we can observe from~\autoref{fig:stencil_size} that using higher PHS orders $k$ degrades stability, as we will further discuss in~\autoref{sec:parametrs}.

\begin{figure}[h!]
	\centering
	\includegraphics[width=1\linewidth]{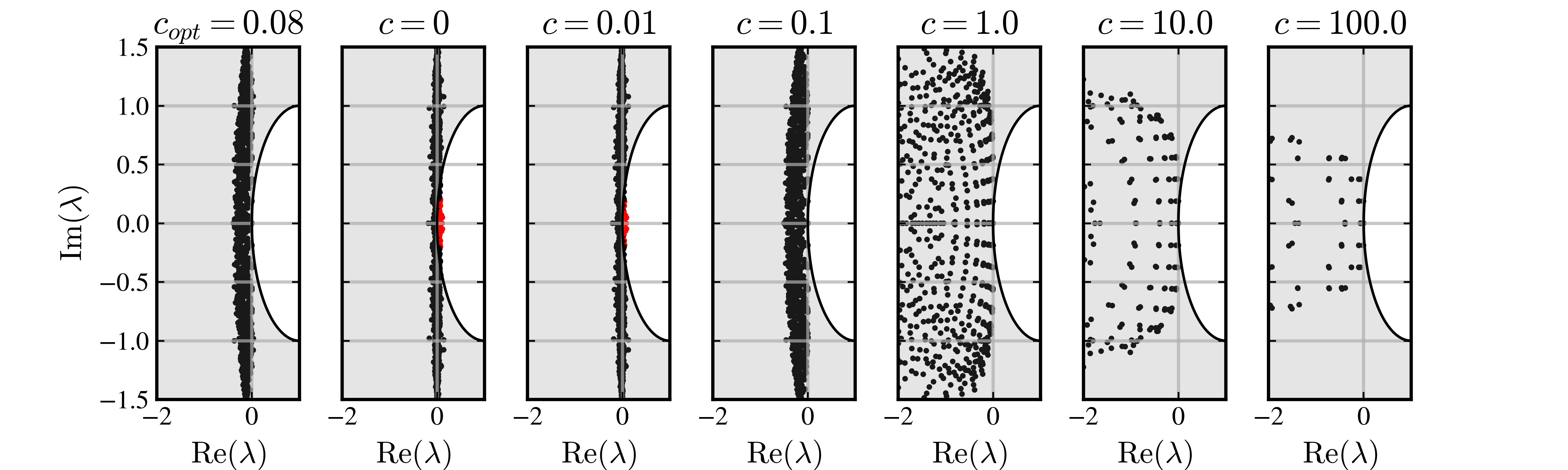}
	\caption{
		Eigenvalue spectra of the stabilised advection operator $\hat{\mt{D}}_h$~\eqref{eq:advection_op}	
		using $N \approx 10^3$ with respect to different values of $c$ for $\alpha = 2$. The shaded region in the figure is the implicit Euler stability region. The monomial order for approximation of $\hat{\mt{D}}_h$ operator is set to $m=2$ and $\Delta t = 10^{-3}$. The black dots represent stable eigenvalues, whereas the red dots represent eigenvalues that lie outside of the stability region. Eigenvalues are scaled with $h$ for ease of visualisation.
}
	\label{fig:eigenvalues_spectra}
\end{figure}

With stencil size proving unable to fully stabilise the system we move onto hyperviscosity. First we examine the eigenvalues of advection operator for a range of hyperviscosity scaling constant $c$. We expect the eigenvalues to progressively shift into the stable region as we increase $c$. This is demonstrated in~\autoref{fig:eigenvalues_spectra}, where eigenvalue spectra of the advection operator $\hat{\mt{D}}_h$~\eqref{eq:advection_op} using $N \approx 10^3$. As expected, we observe that as $c$ increases, more eigenvalues shift into the stable region, where $c_{opt}$ stands for the optimal shift. This optimal shift is defined in \autoref{section:c} as the minimal $c$, where all eigenvalues lie in the stability region and will later be determined by the Algorithm~\autoref{alg:bisection}. However, we also note that the spectra become more scattered, indicating the dissipative properties of the hyperviscosity operator. Consequently, it can be anticipated that as $c$ becomes larger, more physical properties of the system are lost.

\begin{figure}[h!]
	\centering
	\includegraphics[width=1\linewidth]{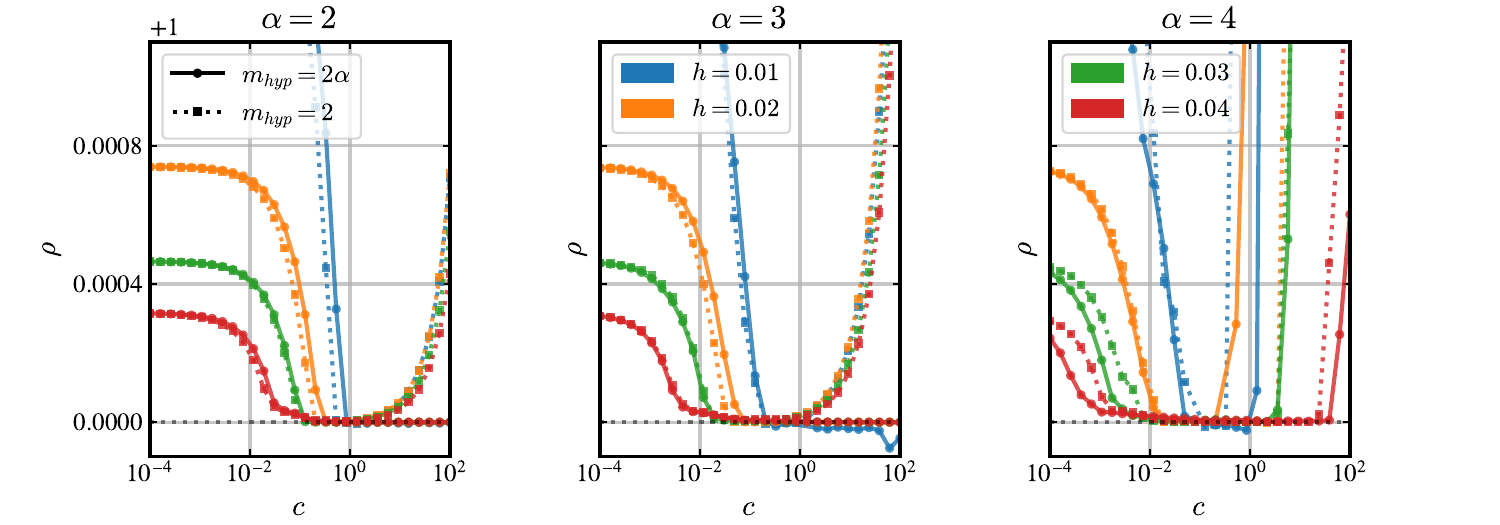}
	\caption{Spectral radius $\rho$ of the evolution matrix $\mt{G}_h$~\eqref{eq:dynamical_system_hv} at $\Delta t = 10^{-4}$ with respect to $c \in [10^{-4}, 10^2]$ for internodal distances $h \in \{0.01, 0.02, 0.03, 0.04\}$ and orders of hyperviscosity $\alpha \in \{2,3,4\}$. The solid line shows $\rho$ calculated with the suggested monomial order for the
    RBF-FD approximation of the hyperviscosity operator while the dotted line displays effects of using a reduced monomial order of $m=2$ with stencil size $n = 30$.}
	\label{fig:advection_spectral_radius}
\end{figure}

Finally, we take a more detailed look at the relationship between the spectral radius $\rho(\mt{G}_h(\gamma(c)))$ and $c$ shown in \autoref{fig:advection_spectral_radius} for different hyperviscosity orders $\alpha$ and internodal distances $h$. At this point we focus on the solid lines, while the dotted lines show the reduced monomial order results discussed in \autoref{sec:parametrs}. The stabilising effect of hyperviscosity appears to be well behaved for the pure advection case i.e. once $c$ is large enough to remove spurious eigenvalues the system remains stable with $\rho \lesssim 1$. Problems only occur in the $\alpha = 4$ case with $\rho$ diverging as we increase $c$. The fact that this is exacerbated by decreasing $h$ indicates that the issue might be caused by increasingly computationally ill-conditioned approximation of the high order derivatives. Decreasing $h$ also increases the required $c$ to reach $\rho \lesssim 1$, which is somewhat unexpected since $h$ scaling is already included in $\gamma(c)$. Further examination is required to determine whether this is caused by the stabilising effect of additional numerical diffusion inherent to larger internodal distances and whether additional scaling for $\gamma$ would be beneficial.

\subsubsection{How to determine the parameter $c$ in practice}
\label{section:c}
The system is stable if $\rho (\mt{G}_h(\gamma(c))) \le 1$. We aim to introduce the least amount of numerical diffusion while maintaining stability. This implies that we are searching for the lowest $c$ (denoted as $c_{opt}$) that produces a stable system
\begin{align}
	c_{opt} = \min \{c \geq 0 \,: \, \rho(\mt{G}_h(\gamma(c))) \le 1 \}.
	\label{eq:set_gamma}
\end{align}

To determine $c_{opt}$, one must solve equation~\eqref{eq:set_gamma}. The simplest approach is to use the bisection method in log-scale, as detailed in Algorithm~\autoref{alg:bisection}. Log-scale is used to improve convergence properties when initial bounds are multiple orders of magnitude apart and to keep the stopping criterion $\Delta$ consistent across the entire range. We use $\Delta=10^{-2}$ resulting in $\sim1\%$ accuracy for the resulting $c_{opt}$. The lower bound is selected in the unstable regime so that $\rho (\mt{G}_h(\gamma(c_{lo}))) > 1$ and the upper in the stable, satisfying $\rho (\mt{G}_h(\gamma(c_{hi}))) \le 1$. The initial bounds were determined with a coarse sweep around an initial guess, but the selection is not trivial for cases with narrow stability regions. Further work is required to devise a more optimal method.

The most computationally intensive part of the algorithm is the stability criterion on line~\autoref{alg:stability_crit} where eigendecomposition of the evolution matrix $\mt{G}_h(\gamma(c_{lo}))$\footnote{The matrix can be decomposed as $\mt{G} = \mt{S} + \gamma($c$) \mt{H}$ to avoid reconstructing the entire matrix whenever $c$ changes.} is required to compute the spectral radius $\rho$. Since we are only interested in computing the eigenvalue with the largest magnitude, the Implicitly Restarted Arnoldi Method (IRAM) algorithm~\cite{Sorensen1997} is employed to improve efficiency but the computational cost of evaluating $\rho$ remains high. Although the computational complexity of IRAM depends on matrix properties and algorithm parameters it can be roughly estimated as $\mathcal{O}(k^2N)$, where $k$ is the number of calculated eigenvalues. We use the iterative algorithm to calculate the $k=40$ largest eigenvalues by magnitude to be reasonably certain that the maximal value was captured since the algorithm does not guarantee convergence to the extremal (largest or smallest) eigenvalue. Note that an alternative stability criterion could be used on line~\autoref{alg:stability_crit}, e.g. a weaker criterion where an energy norm is evaluated to determine stability after propagating the system for a set time~\cite{pozun2025}.

\begin{algorithm}
	\scriptsize
	\caption[Bisection algorithm for $c_{opt}$]{Bisection algorithm for $c_{opt}$}
	\label{alg:bisection}
	\vspace{2pt}
	\textbf{Input:} A method that returns evolution matrix with $c$-dependent hyperviscosity $\mt{G} : \R \to \R^{N \times N}$. \label{alg:evoMatrix}\\
	\textbf{Input:} A method that determines the spectral radius of a matrix $\rho : \R^{N \times N} \to \R$. \\
    \textbf{Input:} Initial guesses for bounds with $c_{lo}$ in the unstable and $c_{hi}$ in the stable regime. \\
    \textbf{Input:} Stopping criterion $\Delta$. \\
	\textbf{Output:} Estimation for $c_{opt}$. \\
    \begin{algorithmic}[1]
		\Function{stabilisationBisection}{$\mt{G}$, $\rho$, $c_{lo}$, $c_{hi}$, $\Delta$}
		\While{$\ln(c_{hi}) - \ln(c_{lo}) > \Delta$}
        \State $c_{mid} \gets \exp((\ln(c_{lo}) + \ln(c_{hi}))/2)$
		\If{$\rho(\mt{G}(c_{mid})) > 1$} \label{alg:stability_crit}
		\State $c_{lo} \gets c_{mid}$
		\Else
		\State $c_{hi} \gets c_{mid}$
		\EndIf
		\EndWhile
		\State \Return $c_{hi}$
		\EndFunction
	\end{algorithmic}
\end{algorithm}

\subsection{Consistency and order of hyperviscosity}
\label{sec:consistency}

Next, we would like to investigate the consistency of the scheme. That is, our approximation of the solution must converge to the correct solution. For simplicity, again assume that $\mathcal{N}$ is a linear operator. Let $u$ be sufficiently regular and solve the posed problem, that is 
\begin{align}
    \partial_t u(\vec x,t) +\mathcal{N}(u(\vec x,t)) = 0.
\end{align}
Then it also holds for many classical PDEs that as $h \to 0$ the family of solutions $\{ u_\gamma \}$ satisfying \eqref{eq:hipscheme} converges strongly $u_\gamma \to u$
since $\gamma = ch^{2\alpha}$. This result extends the classical result (for artificial viscosity), see for example \cite{Tadmor2004,Hussein2020}. Furthermore, it ensures the consistency of stabilisation in a continuous sense. To preserve the order of approximation of the discrete scheme, the stabilisation term must be consistent and vanish faster than the expected convergence rate. Assume that our base convergence rate is $\mathcal{O}(h^p)$, then the hyperviscosity term must vanish with at least the same order to preserve the order of accuracy. To study consistency, we introduce the local interpolation error estimate on a Voronoi region $K_i$ around $\vec x_i \in X$. Fix $t \in [0, T]$ and define $u(\cdot) = u(\cdot, t)$. We define an interpolation operator that maps our functions to the trial space $V_h =  \mathrm{span} \{ \Psi_1, \dots, \Psi_N \}$ by ${I}_h: C(\Omega) \to V_h$ and $ I_h u = \sum_{i=1}^N \Psi_i u(\vec x_i)$. Consider a polyharmonic spline approximation of order $k > 2 \alpha$, augmented with monomials of degree $m$ as explained in \autoref{sec:rbffd}. Provided that the interpolant exists, we have the following error estimate for sufficiently small $h$ and $u \in W^{m+1}_\infty(\Omega) \cap C^{2\alpha}(\Omega)$:
\begin{equation}
\label{eq:errorEst}
\begin{aligned}
    \norm{ \Delta^{\alpha}(I_h u - u)}_{L^\infty(K_i)} &\le C_i h^{m+1 -2\alpha} \abs{u}_{W^{m+1}_\infty(K_i)}. \\
\end{aligned}
\end{equation}
For the $m \geq 2 \alpha$ case, this is a well-known result~\cite{Tominec2021}. If $u \in W^{2 \alpha+1}_\infty(\Omega) \cap C^{2 \alpha}(\Omega)$, then the above also holds for $m < 2\alpha$ as proven in the~\ref{sec:appendixErrorProof}.

Since we would like to study consistency in the $\norm{\cdot}_h$ norm, we must relate it to the $\norm{\cdot}_{L^\infty(\Omega)}$ norm. This is simple for $w: \Omega \to \R$, for which the essential supremum is equal to the supremum:
\begin{equation}
    \label{eq:normConversion}
    \norm{w}^2_h = \frac{|\Omega|}{N} \sum_{i=1}^N w(x_i)^2 \le \frac{|\Omega|}{N} N \max_{1 \le i \le N} w(x_i)^2 = |\Omega| \norm{ w}_\infty^2  \le C \norm{w}^2_{L^\infty(\Omega)}.
\end{equation}
Our local estimates $w = \Delta^\alpha(I_h u - u)$ are only defined on $\overline{K_i}$ but the measure of Voronoi region boundaries is $0$ and therefore $\norm{w}_{L^\infty(\Omega)} = \max_{1 \le i \le N} \norm{w}_{L^\infty\left(\overline{K_i}\right)}$.
Since $w$ is continuous, we have pointwise control via $L^\infty$ norm inside each Voronoi region $K_i$. Hence we again see that $\norm{w}_h \le C \| w\|_{L^\infty(\Omega)}$ holds.
This now allows us to estimate the norm of the stabilisation term at time $t$.
Let $u$ solve~\eqref{eq:hipscheme} and $u_h \in V_h$ be its stable numerical approximation of order $h^p$, furthermore, assume that $m \ge p-1$. We infer by triangle inequality that
\begin{equation}
\begin{aligned}
    \|{\gamma \Delta^{\alpha} u_h}\|_h &= \gamma \norm{\Delta^{\alpha} u_h  + \Delta^{\alpha} u - \Delta^{\alpha} u +  \Delta^{\alpha} I_h u - \Delta^{\alpha} I_h u}_h \\
    &\le ch^{2\alpha} (\norm{ \Delta^{\alpha} u_h - \Delta^{\alpha} I_h u}_h + \norm{ \Delta^{\alpha} I_h u - \Delta^{\alpha} u}_h + \| \Delta^{\alpha} u \|_h). \\
\end{aligned}
\end{equation}
We can bound the first term as
\begin{align}
    \norm{\Delta^{\alpha} u_h - \Delta^{\alpha} I_h u}_h &= \norm{\Delta^{\alpha}({u_h} - {I_h u})}_h \le Ch^{-2\alpha} \norm{{u_h} - {I_h u}}_h,
\end{align}
where an inverse type-inequality was used~\cite[Lemma 4.5.3]{Brenner2008}, and the second according to inequalities~(\ref{eq:normConversion},~\ref{eq:errorEst}). Furthermore we note that by the stability assumption and timestepping consistency we have $\norm{{u_h} - {I_h u}}_h \le  Ch^p$. This yields
\begin{equation}
\begin{aligned}
    \label{eq:hyperviscScaling}
    \|{\gamma \Delta^{\alpha} u_h}\|_h &\le ch^{2\alpha} (C_1 h^{p - 2\alpha} + C_2 h^{m+1 - 2 \alpha}  + C_3) \\
    &\le (C_1 h^{p} + C_2 h^{m+1}  + C_3 h^{2\alpha}).
\end{aligned}
\end{equation}
Note that the $h^{2\alpha}$ scaling for $\gamma$ cancels out the $ch^{-2\alpha}$ term in operator approximation error estimate~\eqref{eq:errorEst}, thus removing the hard limit for minimal monomial order. For sufficiently large $m$ and $p$, the term $h^{2\alpha}$ dominates, which means that our solution can only preserve accuracy up to $2\alpha$ order.
For the $p$-order accurate scheme, we must select $\alpha \ge p/2$ and $m \ge p-1$.

\begin{figure}[h!]
	\centering
	\includegraphics[width=1\linewidth]{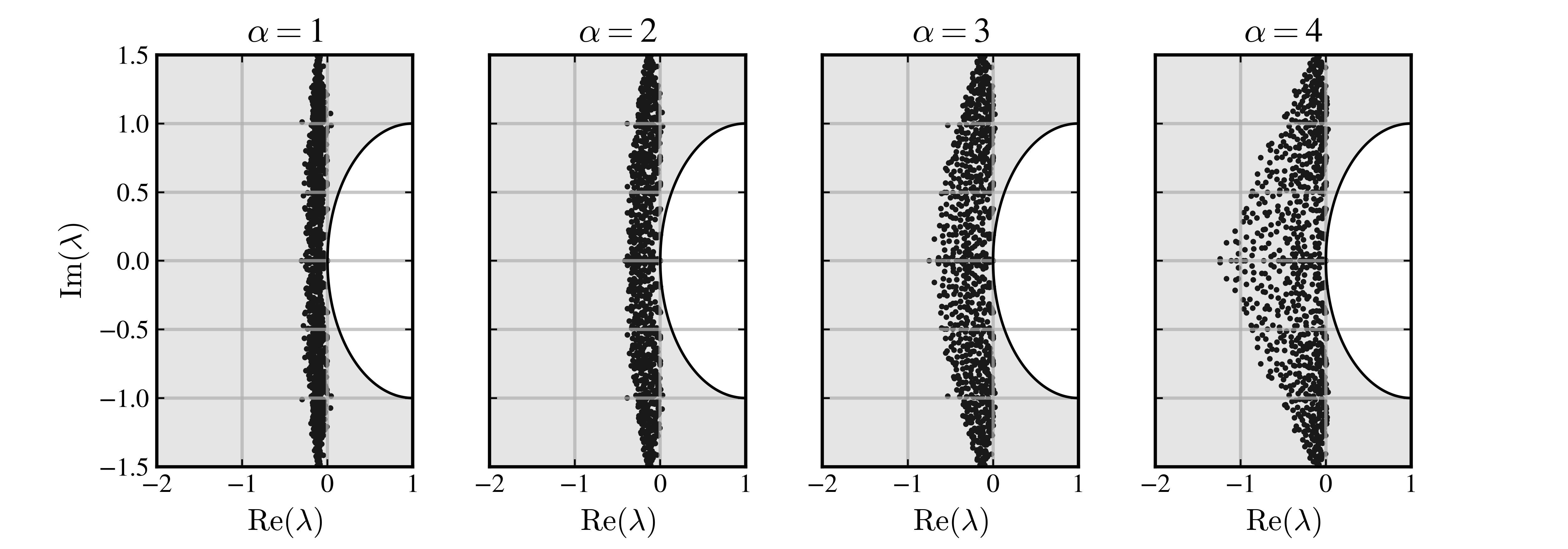}
	\caption{Eigenvalue spectra of the stabilised advection operator $\hat{\mt{D}}_h$~\eqref{eq:advection_op}	
		using $N \approx 10^3$ with respect to different hyperviscosity orders $\alpha$ using $c_{opt}$. The shaded region in the figure is the implicit Euler stability region. Eigenvalues are scaled with $h$ for ease of visualisation.
	}
	\label{fig:eigenvalues_alpha}
\end{figure}

For a second-order approximation, we expect that all orders of hyperviscosity larger than $\alpha \ge 2$ will produce similar results, which we confirm on \autoref{fig:eigenvalues_alpha}. Here, we also note that the approximation of higher-order derivatives becomes increasingly computationally ill-conditioned, since the condition numbers of the differentiation matrix $\mt{D}_h^{hip}$ scale as $h^{-2\alpha}$, hence, lower orders are typically prefered.

\subsection{Parametrisation of the RBF-FD approximation}
\label{sec:parametrs}
One of the problems with hyperviscosity schemes is their large computational complexity. Typically, when approximating operators of order $k$, a monomial order of at least $k$ is required for the first-order operator approximation, i.e., $m = 2\alpha$ is typically used for hyperviscosity. Since the order of the monomial term is inherently connected to the stencil size and subsequently the number of non-zero row entries in the matrix, it has a significant impact on the time complexity of the algorithms for both solving the system and finding the maximal eigenvalues. 

In the previous section, we have shown that PHS alone can be used to control the consistency of the hyperviscosity operator, assuming that the local interpolation systems are unisolvable. This suggests that we can approximate the hyperviscosity operator with much lower monomial orders than otherwise warranted by the high-order derivatives, preserving sparsity and improving computational performance. We also believe that in practice, the bounds that are obtained from inverse inequality~\eqref{eq:invers_inequality} are not strict as argued in \autoref{sec:consistency}. The veracity of this claim and the impact of reduced order approximation for hyperviscosity is first examined in~\autoref{fig:advection_spectral_radius}, where the dotted lines show results with hyperviscosity operators approximated with 2nd order monomials and stencil size $n=30$. In contrast to the full order results, displayed with solid lines, the system now becomes unstable again as $c$ increases. Reduced order results closely match the full order results where $\rho$ first dips below 1 with increasing $c$, which is the region we aim for with the bisection algorithm. Increased stencil size widens the range of $c$ where $\rho \lesssim 1$ and can prevent successful stabilisation if chosen too small. The interplay between the reduced monomial order and stencil size is further analysed in \autoref{sec:linearadvection}.

\begin{figure}[h!]
	\centering
	\includegraphics[width=1\linewidth]{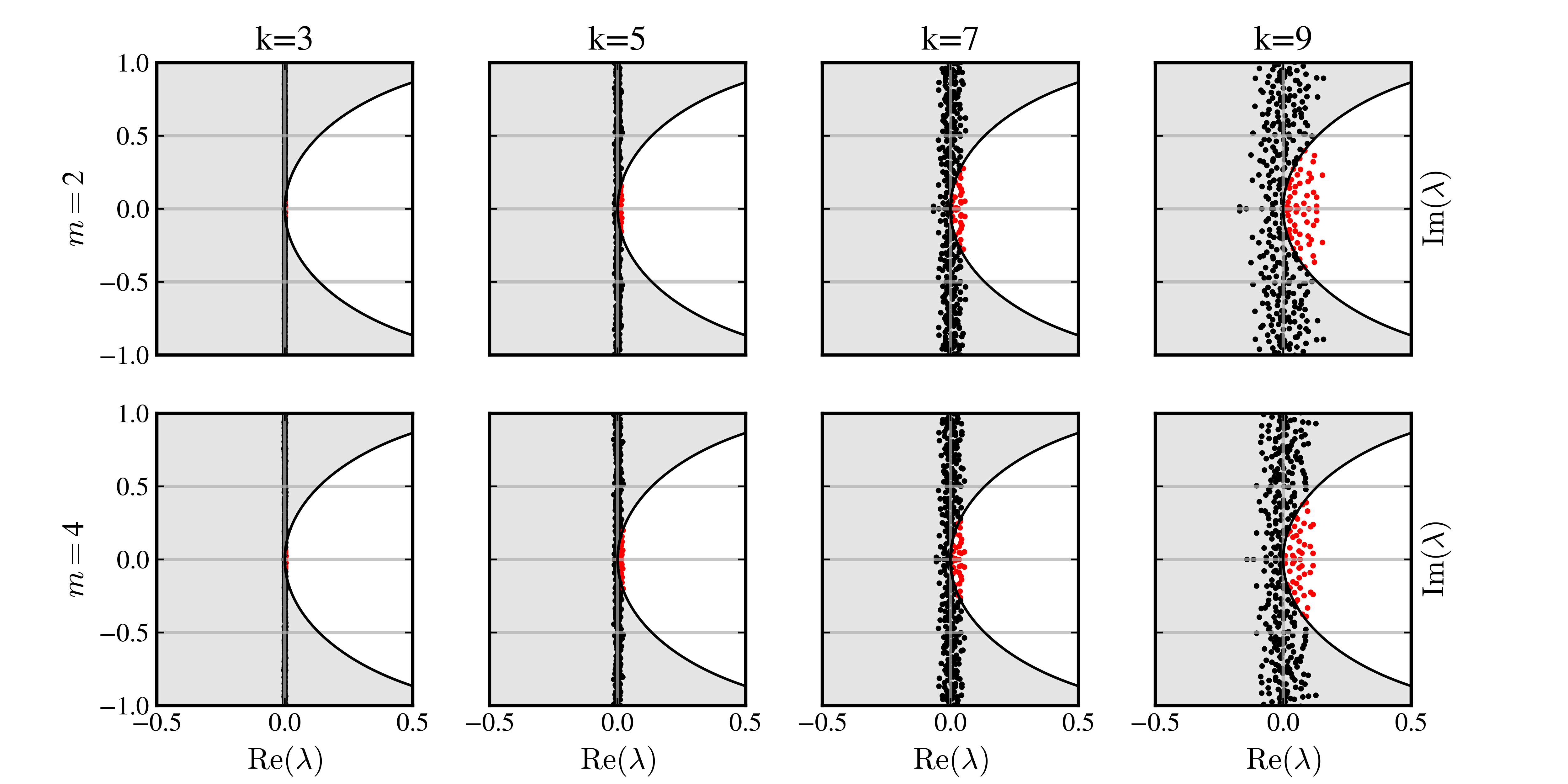}
	\caption{Eigenvalue spectra of the advection operator $\mt{D}_h$ for different orders of polyharmonic splines $k$ at $h=0.03$. The eigenvalues are scaled by $\lambda = h\tilde \lambda$. The shaded region represents the implicit Euler stability region.}
	\label{fig:rbf-fd-spectra}
\end{figure}

It is also essential to select the correct polyharmonic spline order, $k$, to ensure the desired stability properties of the semi-discrete system. The authors of~\cite{shankar2019} proposed a formula linking the monomial order $m$ of the approximation to the polyharmonic spline order via $k = 2m + 1$. This choice is justified by the fact that the system~\eqref{eq:local_interp_system} is guaranteed to be solvable when this polyharmonic spline order is used~\cite{Fasshauer2007}.
However, in practice, the system tends to be solvable regardless of this criterion. Specifically, from our experiments, we observed no significant relationship between the monomial order and the required spline order, as demonstrated in~\autoref{fig:rbf-fd-spectra}, where eigenvalue spectra for the advection operator are shown for various combinations of $k$ and augmentation orders.
Moreover, \autoref{fig:rbf-fd-spectra} confirms the findings from \autoref{fig:stencil_size} that increasing $k$ leads to greater instability. This suggests that the lowest value of $k$ that still provides the required regularity should be used, i.e. for the advection operator this corresponds to $k = 3$. Nonetheless, when approximating hyperviscosity, it is crucial to ensure that the local interpolation system possesses sufficient regularity to compute the $2\alpha$-order derivatives. This requirement implies $k = 2\alpha + 1$. 

Based on this discussion we propose using different approximations for constructing the advection and hyperviscosity operators. Reduced order approximation is used for hyperviscosity term since there appears to be no clear downside. Unless otherwise specified we use $k = 3$ and $m = 2$ with the suggested stencil size $n = 12$ for the advection operator, and $k = 2\alpha + 1$ and $m = 2$ with stencil size $n=30$ for the hyperviscosity operator.

\section{Numerical results}
\label{sec:numerical_results}
In this section, we investigate the impact of hyperviscosity on two unstable problems:
\begin{enumerate}
    \item Linear advection equation, where a  $C^\infty(\Omega)$ initial condition is linearly advected under constant velocity field on $\Omega= \mathbb{T}^2$,
    \item Burgers' equation, where $C^\infty (\Omega)$ initial conditions with simple Dirichlet boundary conditions are non-linearly advected.
\end{enumerate}
These test cases have already been extensively used in a stabilisation research setting. A linear advection equation has also been studied in \cite{Tominec2022,shankarhyperviscosity}. However, our results focus on the influence of the $c$ term on the hyperviscosity stabilisation.  Throughout this section, we use the following error norm
\begin{align}
    \norm{e}_{2} := \frac{\norm{\vec{u} - \vec{u_h}}_{2}}{\norm{\vec{u}}_2} ,
\end{align}
where $\vec{u}$ is a vector of analytical solution values sampled at nodes $X$.
To further show the effects of hyperviscosity, we introduce the relative energy
\begin{align}
    \text{Relative energy} := \frac{\norm{u_h}^2_h}{\norm{u}^2_h} .
\end{align}
  The set of computational nodes $X$ is generated using the dedicated meshless node positioning algorithm described in \cite{Slak2019} with a constant density throughout the computational domain parametrised with the internodal distance $h$. 

 \subsection{Linear advection}
\label{sec:linearadvection}
In this section, we consider the classical 2D pure linear advection equation on a periodic domain $\Omega = \mathbb{T}^2([0,1])$ 
\begin{align}
    {\partial_t} u(\vec{x},t) + \vec{\beta}\cdot\nabla u (\vec{x},t) = 0,
\end{align}
where we set $\vec{\beta} = (1, 0)$. To avoid any issues with the regularity of $u$, we consider the following initial conditions
\begin{align}
    u (\vec{x}, 0) = \begin{cases}
        \exp{\left( \frac{1 - R^2(r^2 - R^2)}{r^2 - R^2} \right )},\quad &r < R\\
        0, \quad & r \ge R.
    \end{cases},
\end{align}
where $R = 0.1$ and $r$ is the distance to the centre of the domain $\Omega$. 
In the fully discrete case we consider
\begin{align}
    \vec{u_h}(t_{n+1}) - \vec{u_h}(t_n) = \Delta t\ ( -\mt{D}^{\partial_x}_h \vec{u_h} (t_{n+1}) + (-1)^{\alpha + 1} \gamma \mt{D}_h^{hip}\vec{u_h}(t_{n+1}) )
\end{align}
The regularity of $u(\vec{x}, 0) \in C^\infty (\Omega)$ allows the application of fourth-order hyperviscosity. The monomial order for approximating the derivatives is set to $m=2$. The error norms are measured by comparing the solutions with the initial solution $u(\vec{x}, 0)$, since $u(x,t) = u(x-t,0)$, at time $t = 1, 2, \dots, T$. The objective of this section is to observe how the stabilised RBF-FD method performs on a pure linear advection.

\begin{figure}[h!]
    \centering
    \includegraphics[width=1\linewidth]{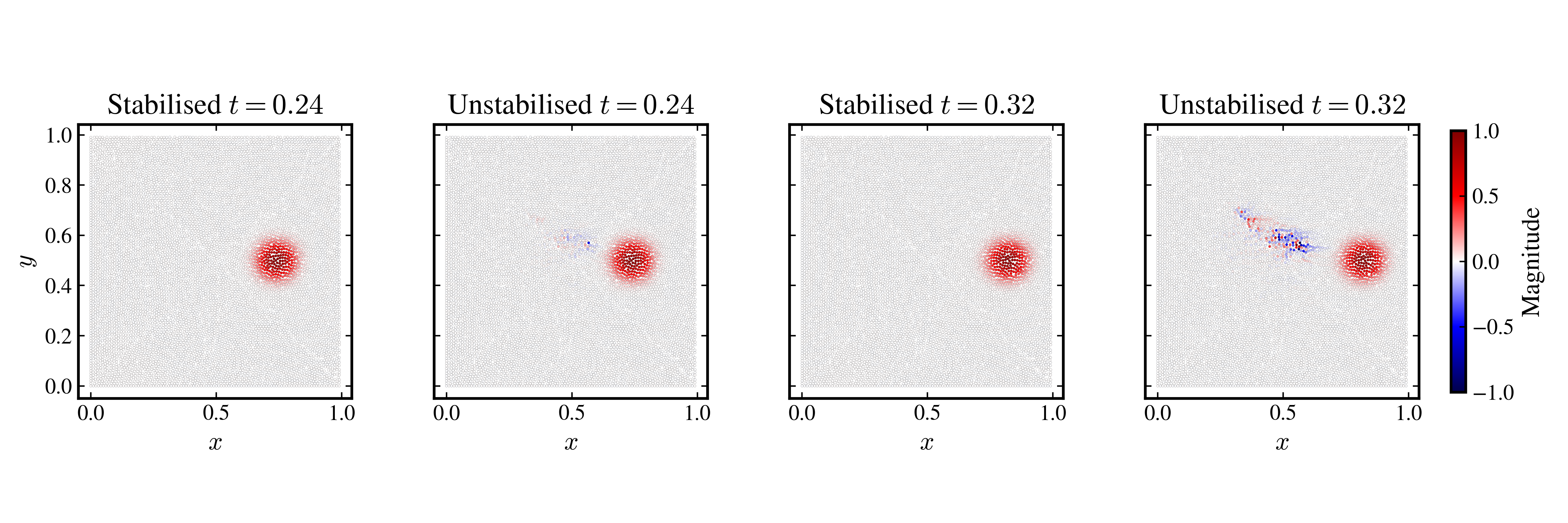}
    \caption{Solution of the linear advection equation at two different times $t \in \{0.24, 0.32\}$ with $h = 0.01$. The stabilised solution (left) admits no visible artefacts, compared to the non-stabilised solution (right) where major instabilities can be observed.}
    \label{fig:advection-comparison}
\end{figure}

We first visually examine 
the solution. In \autoref{fig:advection-comparison}, we compare the hyperviscosity stabilised solution to its non-stabilised counterpart. The stabilised solution uses hyperviscosity of order $\alpha=2$ with $c_{opt}=0.32$ found by the algorithm presented in \autoref{sec:stability_gamma}. 

After a relatively short time $t=0.24$, visible numerical artefacts begin to appear in the non-stabilised solution, whereas no visible artefacts are present in the stabilised solution. Once present, the erroneous field values quickly grow in magnitude, surpassing the advected solution at $t=0.32$ and diverging into invalid numerical values shortly afterwards.

\begin{figure}[h!]
	\centering
	\includegraphics[width=1\linewidth]{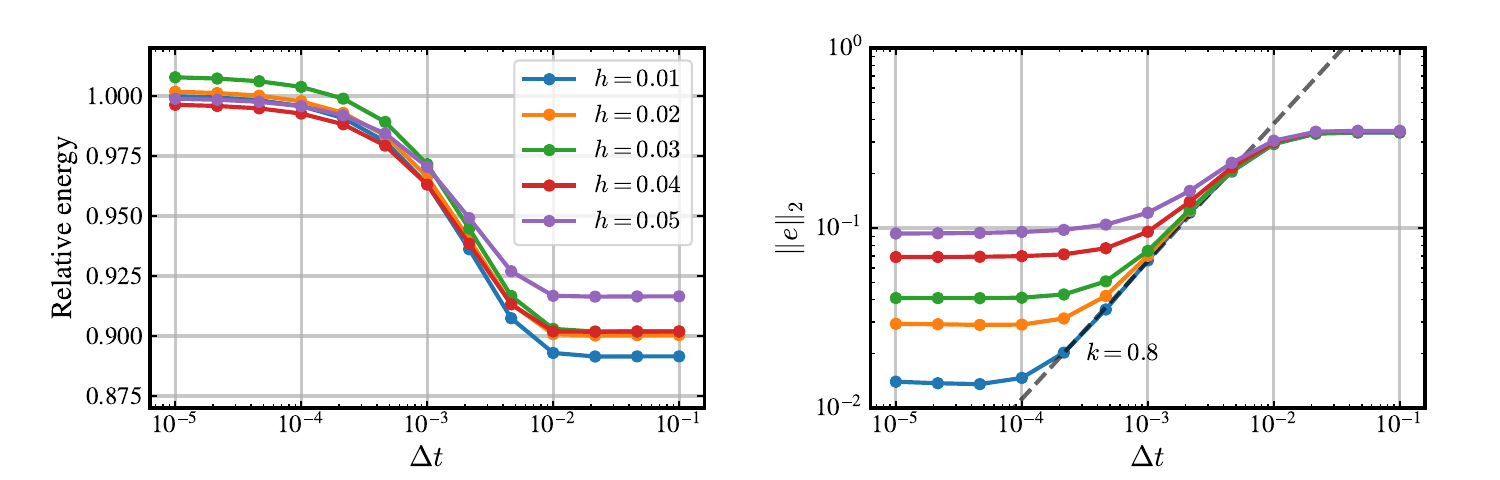}
	\caption{Relative energy and error of numerical simulation of the pure advection with respect to $\Delta t$ at $t = 10$. Stabilised with hyperviscosity of order $\alpha=2$ with a fixed $c = 1$ for clearer comparison. Different lines correspond to different spatial discretisation parameters $h$. The convergence order is fitted based on the middle 5 markers with $\Delta t^k + d$.}
	\label{fig:advection-timestep-convergence}
\end{figure}

Before proceeding, we consider the convergence with respect to timestep $\Delta t$ for different values of $h$ shown in \autoref{fig:advection-timestep-convergence}. The diffusivity caused by the implicit Euler scheme vanishes with $\Delta t \to 0$ and the relative energy, displayed on the left graph, converges towards $1$ indicating energy conservation. This is an encouraging result, since, in general, this is not guaranteed by the stability criterion. On the right graph, showing the relative error, we can observe the convergence behaviour $\mathcal{O}(\Delta t)$. For $\Delta t > 0.01$ the spatial error is dominated by the implicit Euler diffusion, whereas for $\Delta t < 0.01$ we obtain the expected convergence rate of around $1$ before plateauing for a given $h$. Based on this result we use $\Delta t = 10^{-4}$ for following computations unless otherwise specified. 

\begin{figure}[h!]
	\centering
	\includegraphics[width=1\linewidth]{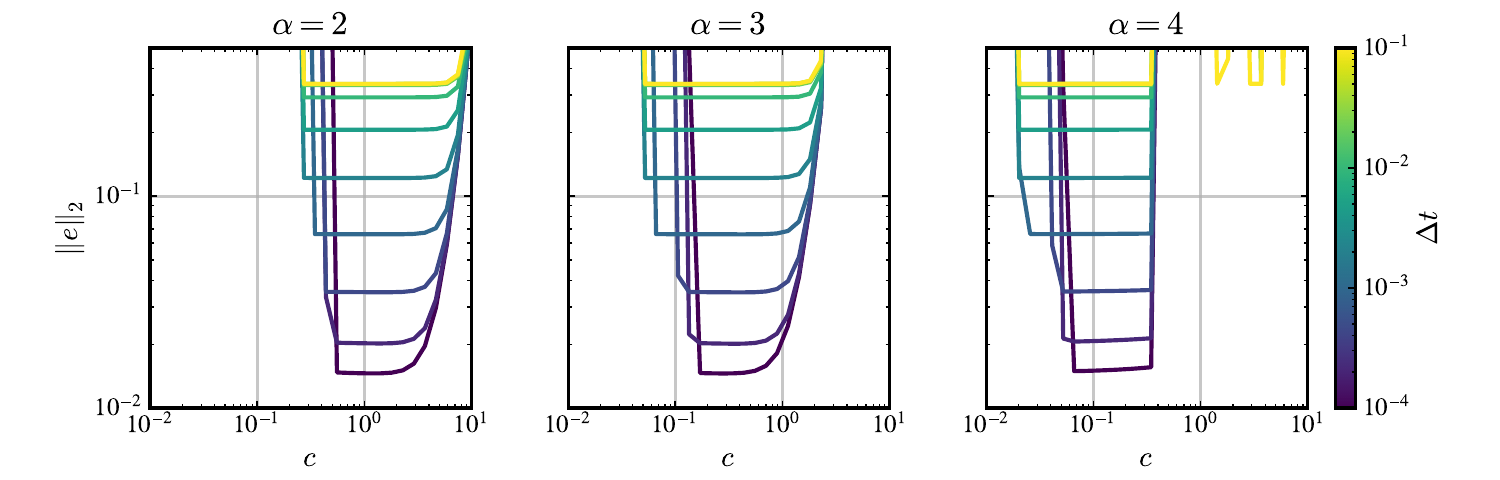}
	\caption{Relative error of pure advection equation at $t=10$ for $h = 0.01$. Different lines correspond to different timestep $\Delta t$ ranging from $0.1$ to $10^{-4}$. Different columns correspond to different orders of hyperviscosity.}
	\label{fig:advection-timestep-c}
\end{figure}

The impact of $\Delta t$ is further scrutinized in \autoref{fig:advection-timestep-c} where we look at its influence on the relationship between the relative error and $c$. The region of stable $c$ shrinks as $\Delta t \to 0$ and less diffusivity is introduced by integration. The relative error is dictated by the 2nd order approximation used for the advection operator and remains similar across different orders of hyperviscosity. This provides us with the first indication that the solution is not disturbed by the stabilising term. 

\begin{figure}[h!]
	\centering
	\includegraphics[width=1\linewidth]{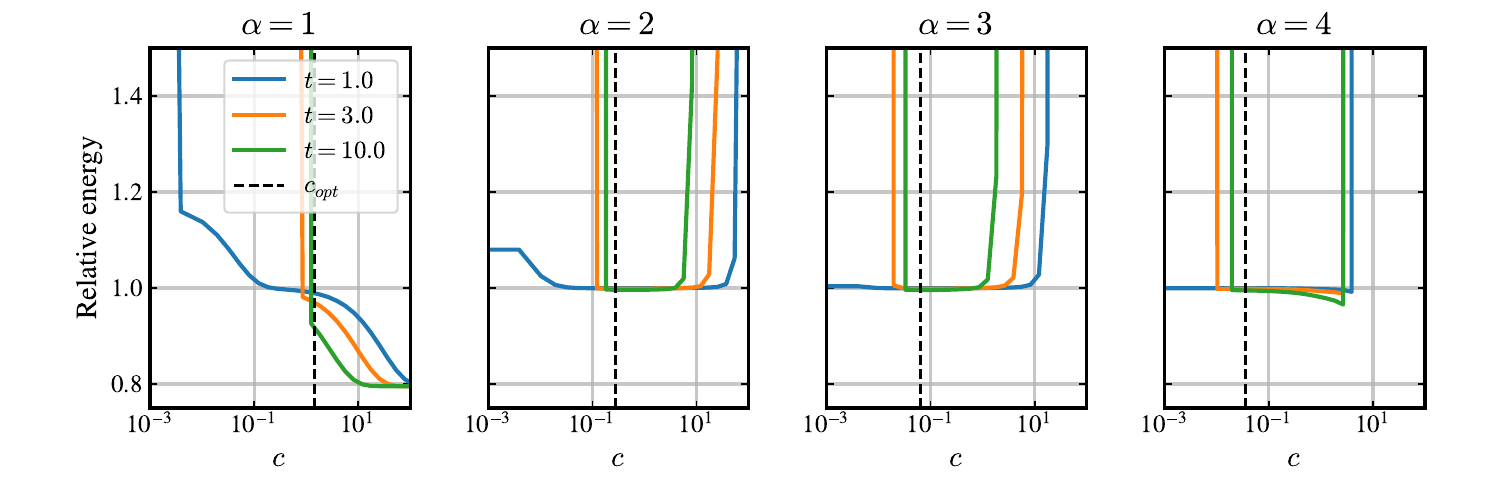}
	\caption{Relative energy of the hyperviscosity stabilized solution discretised with $h=0.02$ with respect to $c$ at $t\in [1, 3, 10]$ for different orders of hyperviscosity.}
	\label{fig:energy-in-time}
\end{figure}

We now discuss the optimal selection of $c$. Although $c_{opt}=0.32$ selected by the algorithm in \autoref{fig:advection-comparison} noticeably stabilises the solution, the optimality of this $c$ remains in question. To shed some light on this issue, we analyse the relative energy of the solution after $1$, $3$, and $10$ rotations, corresponding to
$t \in \{1, 3, 10\}$, for different orders of hyperviscosity $\alpha$. The results of this analysis, performed with $h = 0.02$ ($N \approx 2100$), are presented in~\autoref{fig:energy-in-time}.
First, we observe that for each setup, there is a range of suitable $c$ values where the solution remains stable with minimal dissipation. This region narrows as smaller and smaller instabilities become apparent with passing time and cause solutions to diverge. Most importantly, the $c_{opt}$ (denoted by the vertical dashed line) falls within the stable range in all cases, which is exactly what we hoped for. We note that the $\alpha = 1$ case, analogous to normal diffusion, exhibits significantly stronger drop in relative energy with increasing $c$ when compared to higher orders of Laplacian operator.

\begin{figure}[h!]
	\centering
	\includegraphics[width=1\linewidth]{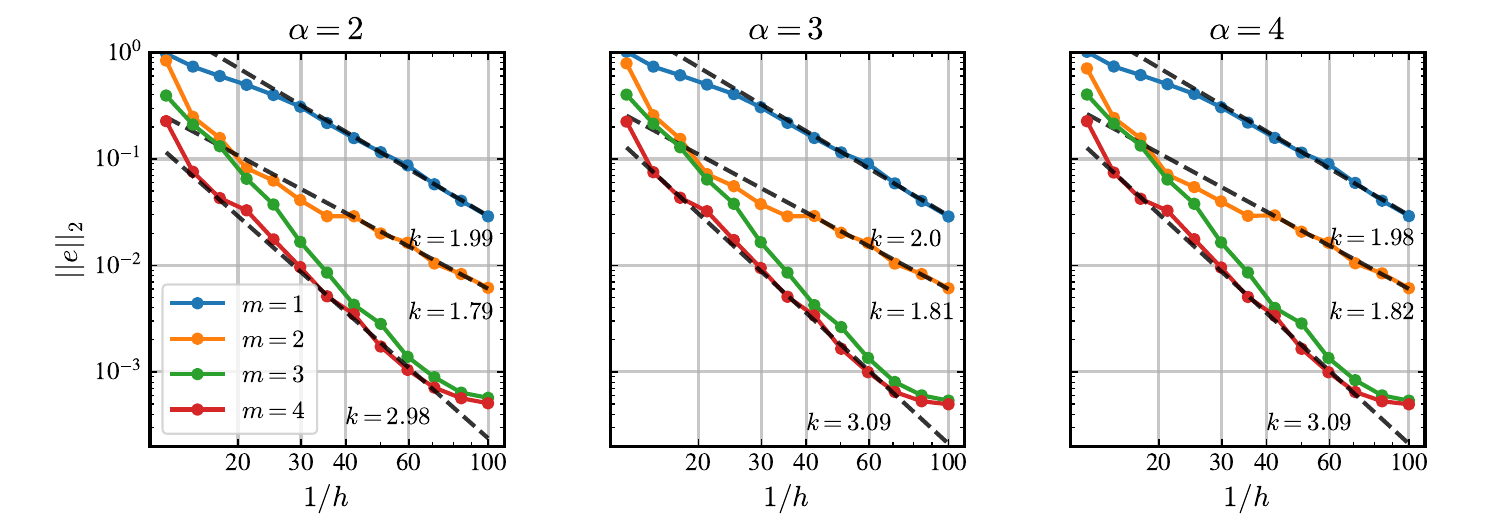}
	\caption{Error convergence under $h$-refinement for linear advection equation at $t=5$ with $\Delta t = 10^{-5}$. Different lines represent different orders of monomial augmentation, while columns represent the order of the hyperviscosity operator. The black lines are fitted based on the ending 8 markers without the last 2 markers with $h^k + d$.}
	\label{fig:advection-consistency}
\end{figure}

To verify that the hyperviscosity operator does not interfere with convergence rates of the scheme, we consider convergence under $h$-refinement for different orders of monomial augmentation $m$ used in approximating the advection term. In \autoref{fig:advection-consistency} we observe that the convergence rates stay relatively stable across different orders of hyperviscosity. In the last two markers for the high order approximations, the effects of timestepping error already dominate over the spatial error.
For $m \in \{2, 4\}$, we obtain 1 order of convergence less than is typically expected. The observed effect is similar to \cite{KolarPoun2024}. The overall convergence rates agree with the observed convergence rates in the RBF-FD setting, most importantly, the diffusion caused by hyperviscosity does not affect the convergence rate.

\begin{figure}[h!]
	\centering
	\includegraphics[width=1\linewidth]{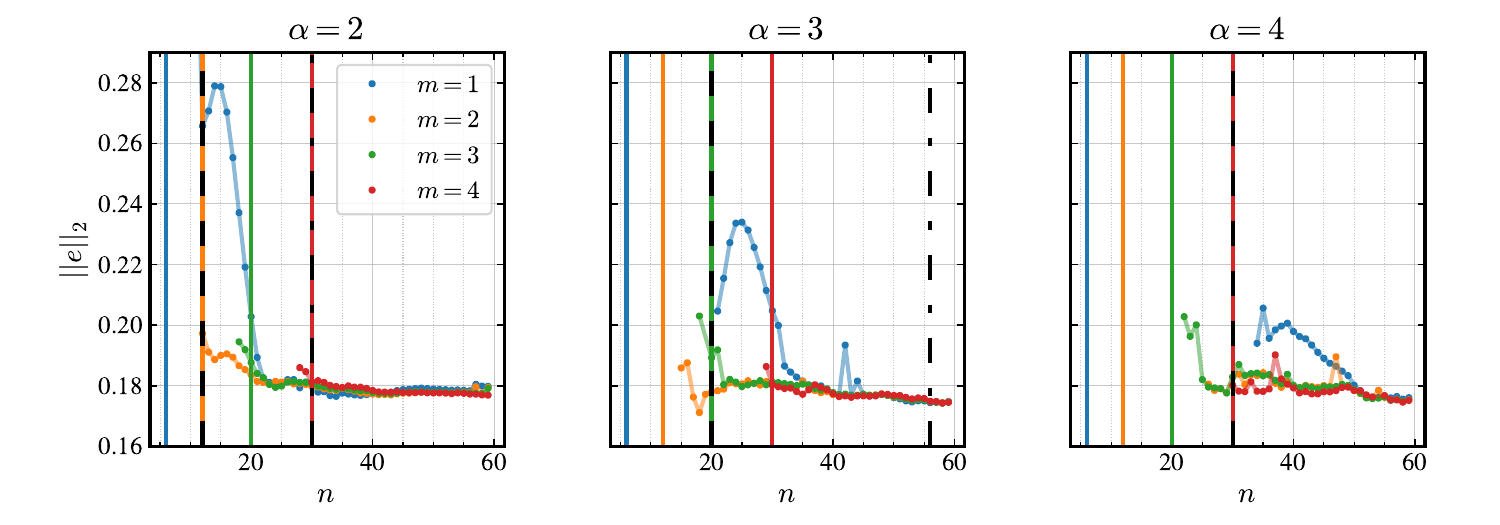}
	\caption{Error of pure advection numerical simulation with respect to stencil size of the hyperviscosity operator at $t=10$ and $h=0.02$. Different colours correspond to different orders of hyperviscosity monomial augmentation $m$. Coloured vertical lines display the recommended~\cite{bayona2019} stencil size for that augmentation order. Black vertical lines show the recommended stencil sizes based on $\alpha$. Dash-dotted for the order appropriate $m = 2\alpha$ and dashed for $m = \lceil \frac{k}{2} \rceil - 1$, with $k = 2\alpha + 1$ based on unisolvency bounds. }
	\label{fig:hyperviscosity-stencil}
\end{figure}

Next, we take a look at how the stencil size and different orders of polynomial augmentation for the hyperviscosity operator affect the error of the solution. This is shown in \autoref{fig:hyperviscosity-stencil}.
We notice that the results are in strong agreement with the theory from \autoref{sec:rbffd_approx} and \autoref{sec:consistency}. We are able to successfully stabilise advection with augmentation orders that are much lower than the approximated derivative orders as indicated by the hyperviscosity consistency estimate~\eqref{eq:hyperviscScaling}. This comes at the cost of increased stencil size though. Stable solutions require stencil size that is roughly in line with the suggested stencil size for $m = \lceil \frac{k}{2} \rceil - 1$ required to prove unisolvency for PHS with $k = 2\alpha + 1$. Even with this increased stencil size requirement there are still significant computational benefits. In $\alpha = 4$ case we require $n \approx 30$ instead of $n \approx 90$ that is suggested for $m = 8$.
The error with $m=1$ tends to dominate the overall error, which is also consistent with the theory, since for that case $m=p-1$ and the $\mathcal{O}(h^{m+1})$ term in our hyperviscosity consistency estimate becomes relevant.

\begin{figure}[h!]
	\centering
	\includegraphics[width=1\linewidth]{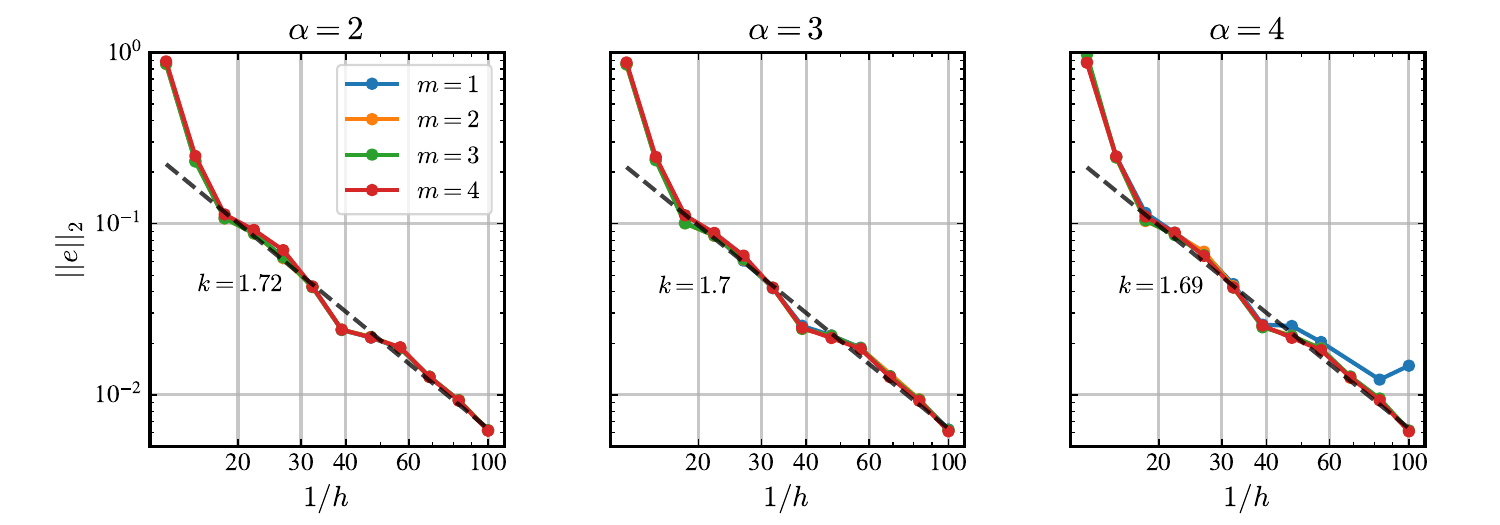}
	\caption{Error of pure advection numerical simulation with respect to $1/h$ at $t=5$. Different lines correspond to different orders of hyperviscosity monomial augmentation $m$.}
	\label{fig:hyperviscosity-consistency}
\end{figure}

At last, we take a look at the consistency of the hyperviscosity operator by varying the monomial order in the approximation of the hyperviscosity operator. According to \autoref{fig:hyperviscosity-stencil}, we use stencil size $n = 35$ to be able to fairly compare all of the considered orders. We set $\Delta t = 10^{-5}$ and perform a similar experiment to the above. In \autoref{fig:hyperviscosity-consistency} we can see that as $h \to 0$ the errors for different approximations (monomial orders) of the hyperviscosity operator stay relatively similar, i.e. the error of approximation stays consistent regardless of the order of monomial augmentation. The only apparent outlier is the severely reduced order $m=1$ approximation for $\alpha = 4$.

\subsection{Non-linear advection: Burgers' equation}
\label{sec:burgers}
In this section, we consider a 2D non-linear advection-diffusion equation, the viscid Burgers' equation 
\begin{align}
    \partial_t \vec{u} (\vec{x},t) + \vec{u} (\vec{x},t) \cdot \nabla \vec{u}(\vec{x},t) = \re^{-1} \Delta \vec{u}(\vec{x},t),
\end{align}
where $\vec{u}$, $\re$ stand for velocity and the Reynold's number, respectively. The problem is considered on a unit square $\Omega = [0,1] \times [0,1]$ with Dirichlet boundary conditions taken from closed form solution 
$\vec f: \Omega \times [0, T] \to \R^2$ from~\cite{Fletcher1983,Zhu2010}, 
\begin{align}
	\vec{f}(\vec{x}, t) =
	\begin{bmatrix}
		\frac{3}{4} - \frac{1}{({4(1+\exp{\frac{\re}{32} (-t-4x_1+4x_2)}))}} \\
		\frac{3}{4} + \frac{1}{({4(1+\exp{\frac{\re}{32} (-t-4x_1+4x_2)}))}}
	\end{bmatrix}.
\end{align}
Again, we solve the problem using implicit time stepping, where the problem is numerically linearised by using the previous time step velocity in the advection term
 \begin{equation}
 \begin{aligned}
 	{\vec{u_h}(t_{n+1}) - \vec{u_h}( t_{n})} =& \Delta t(-\vec{u_h} ( t_{n}) \cdot \mt{D}_h^\nabla \vec{u_h}(t_{n+1}) + \re^{-1} \mt{D}_h^\Delta \vec{u_h}( t_{n+1}) \\& + (-1)^{\alpha +1} \gamma  \mt{D}_h^{hip}\vec{u_h}(t_{n+1})),
 \end{aligned}
 \end{equation}
The linearisation of the advection term causes the matrix $\mt{G}_h$ to become time-dependent; i.e., $c$ is also time-dependent, and it should therefore be recomputed multiple times during the simulation. We also note that this is not a linear Cauchy problem, and the Lax-stability is not a correct notion of stability for this problem. However, we assume that boundary conditions and nonlinearity do not play a crucial role and still use it as a valid criterion. The equation becomes more unstable with increasing $\re$, as advection begins to dominate over diffusion. Moreover, the solution exhibits increasingly pronounced shock-like fronts with thin boundary layers near the shocks, where the Gibbs phenomenon may occur. Here, we aim to demonstrate that in the moderately diffusive regime, a stable solution can be obtained using hyperviscous RBF-FD, whereas the non-stabilised RBF-FD diverges. 

We continue to use $m=2$ to approximate the equation operators as there is little to be gained with higher order approximation of the relatively sparsely discretised narrow front~\cite{Jancic_Kosec_2024} that appears at higher $\re$.

In the following discussions, if not specifically stated otherwise, we recompute $c_{opt}$ every $\chi=10$ time-steps.

\begin{figure}[h!]
	\centering
	\includegraphics[width=1\linewidth]{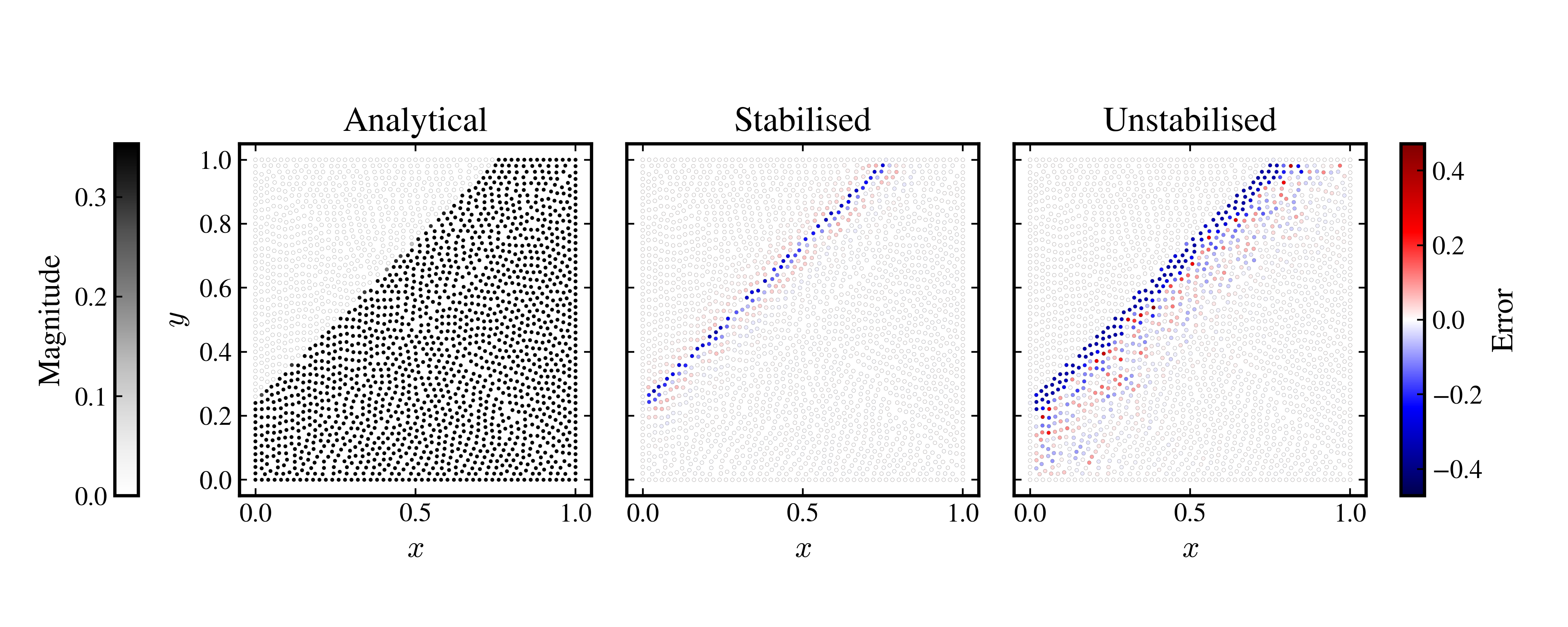}
	\caption{Comparison between stabilised (centre) and non-stabilised (right) Burgers' equation at $t=1$. The error is measured using a pointwise $\ell_1$ norm. The simulation was performed at $\re=5000$ with $h=0.01$ ($N\approx 10^4$) with $\Delta t = 0.01$.}
	\label{fig:burgers-comparison}
\end{figure}

To begin with, we visually compare the stabilised and non-stabilised solutions of the Burgers' equation in~\autoref{fig:burgers-comparison}. The parameters of the scheme are $h = 0.01$ and $\re = 5000$. Although the initial conditions are smooth and the effects of non-linear advection are still relatively small, the non-stabilised solution exhibits considerable numerical artefacts arising from the interplay of the Gibbs phenomenon, advection operator instability, and the inherent characteristics of RBF-FD~\cite{Tominec2025}, whereas these artefacts are less present in the stabilised solution.

\begin{figure}[h!]
	\centering
	\includegraphics[width=1\linewidth]{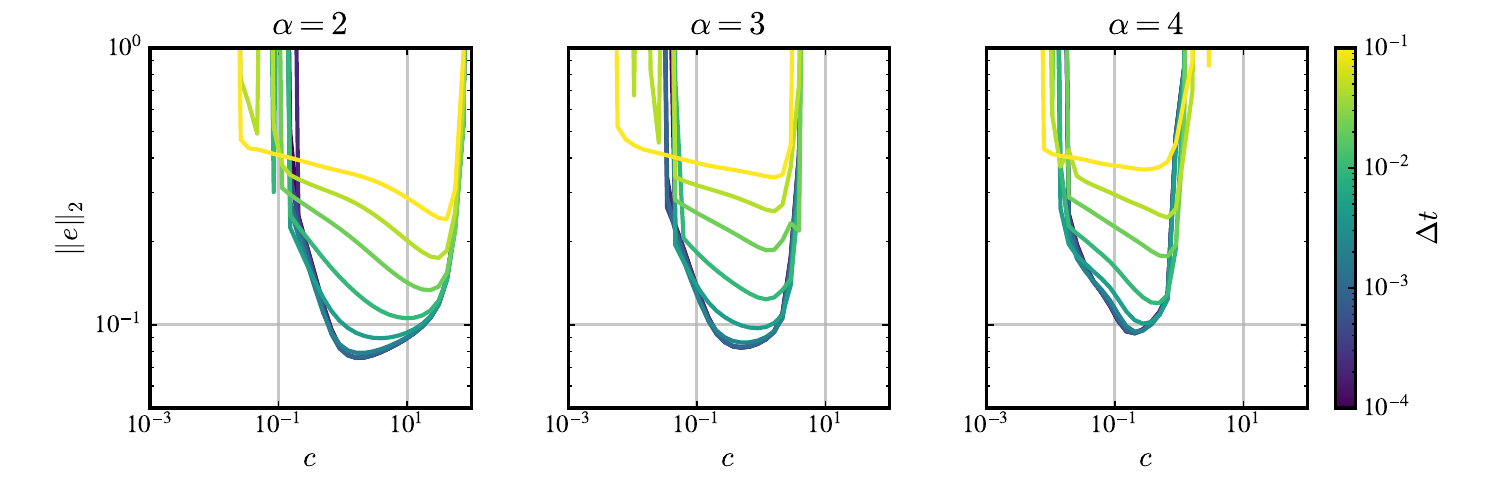}
	\caption{Error for the case of Burgers' equation at $t=1$ for $h = 0.02$. Different lines correspond to different timesteps $\Delta t$ ranging from $0.1$ to $10^{-4}$. Different columns correspond to different orders of hyperviscosity.}
	\label{fig:burgers-dt}
\end{figure}

Akin to the previous section, we begin with the timestep analysis. We again consider the error with respect to $c$ for different $\Delta t$. We choose an unstable case of $\re=10000$, where shocks already form. The results are shown in \autoref{fig:burgers-dt}. This time, the stability area does not shrink as $\Delta t \to 0$, however, it does shrink with the hyperviscosity order $\alpha$. The order $\alpha=2$ results in the lowest error and the largest stability area, which is consistent with the fact that shocks need lower-order smoothing. This prompts us to use the order $\alpha=2$ for the remaining results. We also observe that there exists a $c$ where the error achieves its minimum and this is no longer the lowest stable $c$ that we search for with our algorithm. 

\begin{figure}[h!]
	\centering
	\includegraphics[width=1\linewidth]{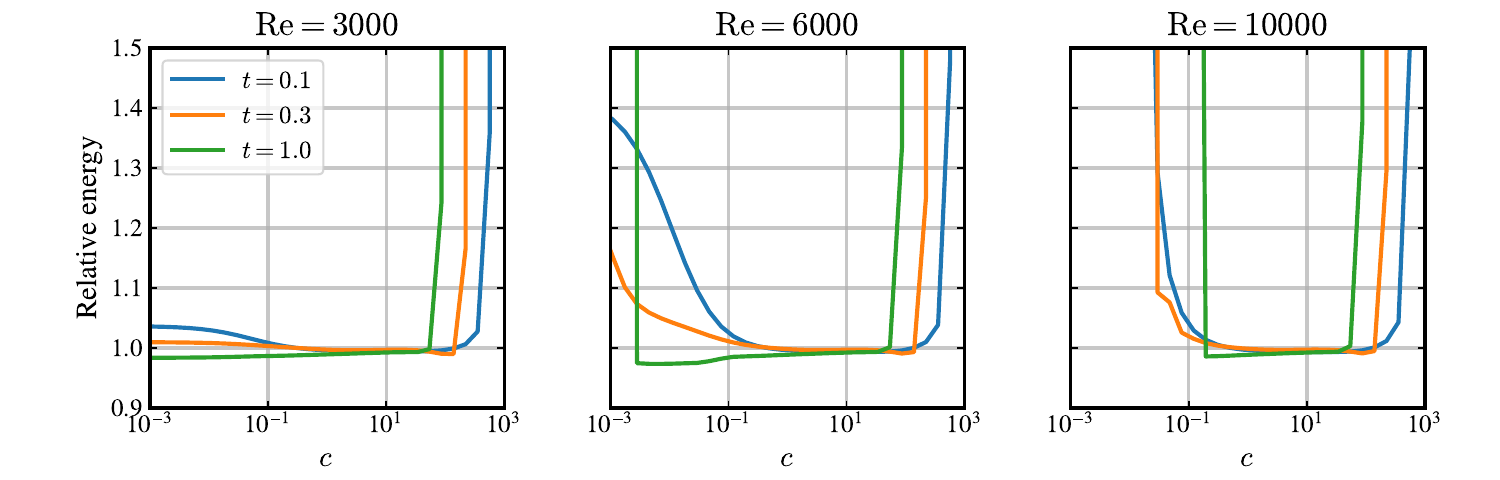}
	\caption{Relative energy with respect to $c$ at $h=0.02$ and $\Delta t=0.01$. Different lines represent sampling at different times of the simulation. The columns show simulation at different Reynolds numbers. }
	\label{fig:burgers-energy}
\end{figure}

 In~\autoref{fig:burgers-energy} we examine the relation between the relative energy and $c$ for different Reynolds numbers at different times of the simulation. Again, we see that there is an interval of suitable $c$ values where the solution remains stable with minimal dissipation. Note that in contrast with the linear advection case shown in~\autoref{fig:energy-in-time} the smallest stable $c$ no longer leads to an optimal result in terms of the relative energy conservation. This loss is caused by insufficient damping of the numerical artefacts that appear when evolving the sharp initial condition. This is most apparent in the intermediate range of $\re$ where the initial condition is sharp enough to cause significant artefacts, while the case is still sufficiently mild to not outright diverge.   

\begin{figure}[h!]
    \centering
    \includegraphics[width=1\linewidth]{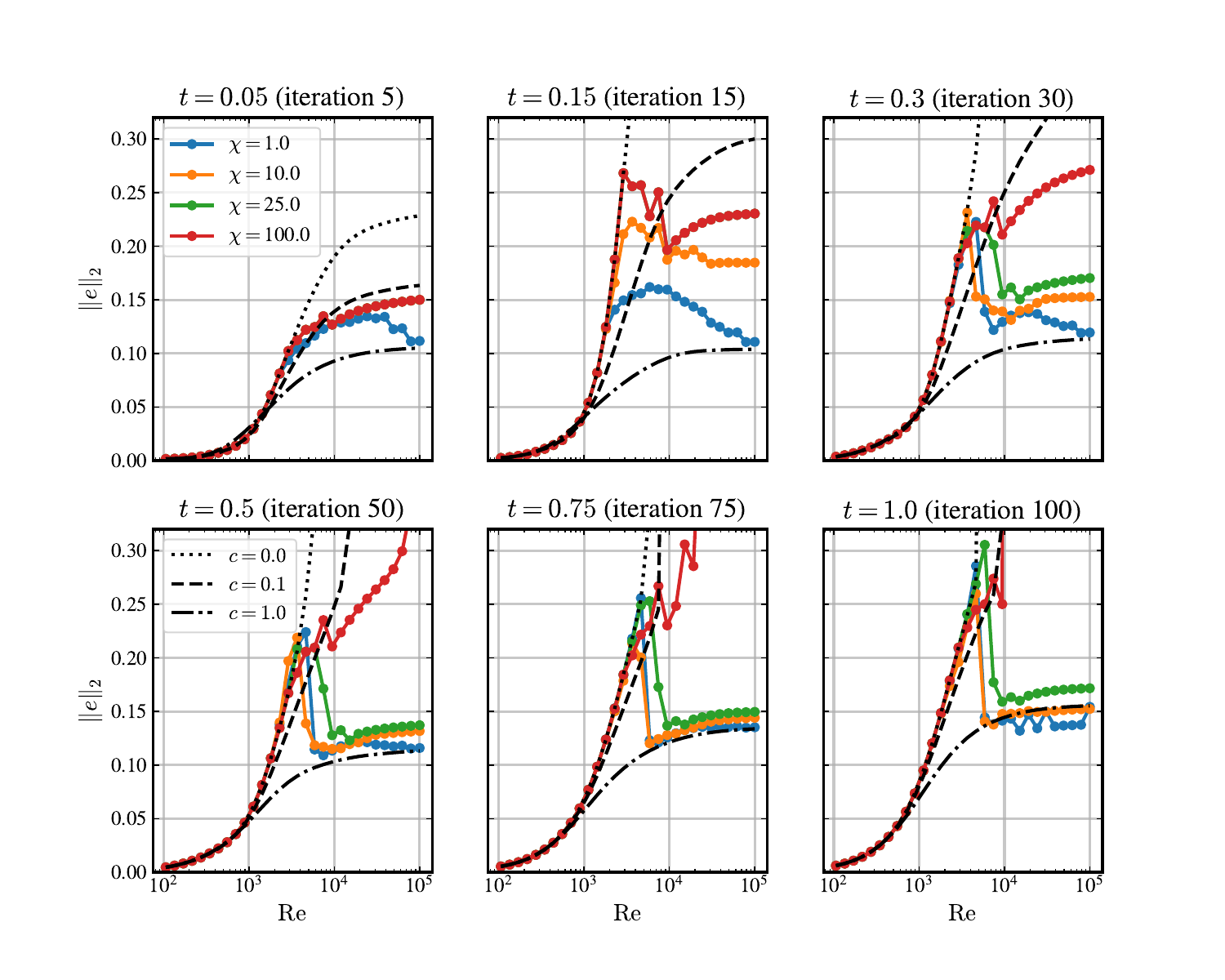}
    \caption{ Time evolution of relative error of Burgers' equation at $\Delta t = 0.01$ and $h=0.02$ with respect to Reynolds number ($\re$). Different colours of lines correspond to $c_{opt}$ at different recomputation frequency $\chi$. The black lines correspond to statically chosen hyperviscosity parameters.}
    \label{fig:burgers-rho}
\end{figure}

As the linearisation of the advection term causes the matrix $\mt{G}_h(t_n)$ to become time-dependent, a natural question arises regarding the required frequency of $c_{opt}$ re-computation. We can assume that the number of times that $c_{opt}$ is recomputed is proportional to the system's dynamics, i.e. the
higher the Reynolds number, the more frequent $c_{opt}$ recomputation is required for optimal stabilisation. This is demonstrated in \autoref{fig:burgers-rho}, where one can see that for Reynolds numbers below 5000, the result quality is weakly coupled with the re-computation frequency and remains similar to the non-stabilised case. This highlights one of the downsides of the algorithm as the cases that do not exhibit exponential growth in their eigenvalues remain un-stabilised even in the regime where a moderate stabilisation displayed with the dash-dotted line would prove beneficial. However, with increasing Reynolds number, the effect of more frequent $c_{opt}$ re-computation becomes evident. Times in the top row of \autoref{fig:burgers-rho} were selected to show the behaviour 5 iterations after the first re-compute for different frequencies while the bottom row displays the rest of the evolution to $t = 1$. We can see that the error of the solution with the newly re-computed $c_{opt}$ clearly drops from the overlap of others that are still stabilised with the original value determined at the first step. Note that the benefits of the more frequent re-computation are still present at $t = 1$ after many lower frequency re-computations indicating that a sufficient re-computation frequency in the initial stages of the evolution is critical.

\begin{figure}[h!]
    \centering
    \includegraphics[width=1\linewidth]{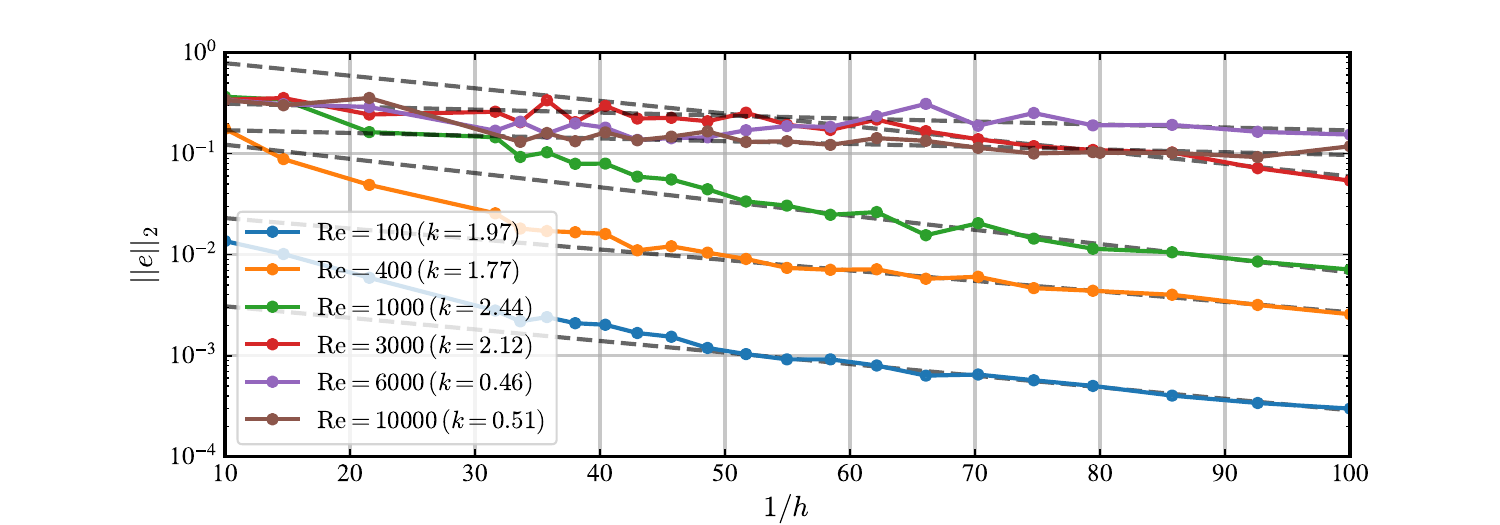}
    \caption{Convergence of the approximation error as a function of $1/h$ for different Reynolds numbers ($\re$) with $\Delta t = 0.01 h$. The error is captured at $t=1$ with $c_{opt}$ recomputed every $10$ iterations. The lines are fitted based on the last 10 markers with $h^k + d$.}
    \label{fig:burgers_consistency}
\end{figure}

Finally, in \autoref{fig:burgers_consistency}, the convergence of the presented solution is shown. We get the expected 2nd order of convergence for low Reynolds numbers $\re \le 3000$ and weak convergence above that is caused by the aforementioned numerical artefacts degrading the increasingly sharp initial wave. 

\section*{Note on implementation}
The numerical procedure is implemented in the Medusa C++ library \cite{slak_medusa_2021}. For the computation of the eigenvalues, we employ the Krylov subspace algorithms \cite{Saad1980} offered by the Spectra C++ library~\cite{spectra} for sparse and the algorithms based on Schur decomposition provided within the function \emph{eigs} in the Eigen C++ library \cite{eigenweb} for dense matrices. 
The code was compiled using g++ (GCC) 11.3.0 for Linux with -O3 -DNDEBUG flags on AMD EPYC 7702 64-Core Processor computer. OpenMP API has been used to run parts of the solution procedure in parallel on shared memory. Post-processing was done using Python 3.10.6.

\section{Final remarks}
\label{sec:discussion}
This paper discusses the hyperviscosity stabilisation scheme in the RBF-FD context from various points of view. The presented analysis was primarily focused on the computation of the adherent hyperviscosity constant $c$ (or $\gamma$) to satisfy the Lax equivalence theorem. In particular, we constructed an algorithm that uses the maximal eigenvalue magnitude of the evolution matrix to iteratively obtain the constant that produces a Lax-stable scheme in the sense that all eigenvalues of the evolution matrix lie inside the stability region. 

Besides the parameter $c$, which effectively governs the magnitude of hyperviscosity, another important parameter, the order of hyperviscosity $\alpha$, determines the vanishing order. As $\alpha$ increases, the order of the derivatives involved increases as  $2\alpha$. Even for low values of $\alpha$, this becomes computationally expensive due to the large stencil requirements. To mitigate this issue, we considered the consistency estimate for the hyperviscosity operator. Specifically, we demonstrated that any monomial order can be used in the RBF-FD approximation of the hyperviscosity operator to produce a consistent scheme. This, in turn, allows for the use of smaller stencil sizes in the RBF-FD approximation. 

In terms of the RBF-FD parametrisation, we showed that increasing the polyharmonic spline order $k$ also increases the instability of the system. On the other hand, when approximating hyperviscosity, sufficient interpolation regularity is required to compute the $2\alpha$-order derivatives, specifically, $k = 2\alpha + 1$. Therefore, we propose using different polyharmonic spline orders for constructing the advection and hyperviscosity operators, namely, $k = 3$ for the advection operator and $k = 2\alpha + 1$ for the hyperviscosity operator.

In the numerical demonstration, we first tackled a linear advection problem. Although the case is relatively simple, its numerical solution exhibits instability, making it an ideal case for studying stabilisation schemes. We initially demonstrated that the proposed scheme provides significantly more stable results compared to the non-stabilised RBF-FD counterpart. We also investigated the interplay between $\alpha$ and $c$, clearly showing that there exists a range of $c$ values for which the solution remains stable while introducing relatively low dissipation. Moreover, the presented algorithm successfully identifies $c_{\text{opt}}$ within this range. The advection case is also used to demonstrate the consistency of the hyperviscosity operator and to verify the theoretical result regarding the use of low monomial orders for high-order hyperviscosity.

To extend the analysis, we addressed a more complex non-linear case, namely the Burgers' problem. First, we repeated similar analyses as in the linear advection case to confirm the adequacy of the proposed solution procedure. However, the main focus of the analysis in this case is the non-linearity of the problem, which makes the parameter $c_{opt}$ time-dependent. We demonstrated how $c_{opt}$ evolves throughout the simulation and to what extent the frequency of its re-computation affects the simulation results.

To sum up, this paper adds another piece to the understanding and development of general hyperviscosity stabilisation. It discusses parameter selection and their influence on the numerical simulation, as well as the core concepts of hyperviscosity itself. Although this paper advances hyperviscosity stabilisation one step closer to practical implementation, several unanswered questions and open topics remain for future work. Starting with the $c$ selection algorithm, more efficient techniques for determining the largest eigenvalues should be investigated, and the bisection algorithm for root finding could be replaced with a more sophisticated approach to reduce the required number of eigenvalue re-computations. There are further opportunities for adaptive stabilisation with locally dependent $c$ and $\alpha$. Although the numerical experiments presented herein are confined to quasi-uniform node distributions, the underlying methodology is theoretically extensible to variable-density discretisations. Further investigation is warranted to evaluate the performance in domains where the nodal spacing $h$ varies spatially.

In time-dependent problems, an adaptive $c$ recomputation could be devised that would follow the dynamics of the solution and recompute $c$ accordingly. For example, in highly unstable cases, it might be practical to periodically compute the eigenvalue with the largest magnitude and then assess whether to reapply the algorithm or not. In future work, we will also consider the recently provided energy method framework~\cite{Tominec2025} to potentially determine $c$ more effectively. 

\section*{Acknowledgements}
\noindent
The authors would like to acknowledge the financial support from the Slovenian Research and Innovation Agency (ARIS) in the framework of the research core funding No. P2-0095 and research projects No. J2-3048 and N2-0275, as well as the Young Researcher programmes PR-12347 and PR-10468.
Funded by National Science Centre, Poland under the OPUS call in the Weave programme 2021/43/I/ST3/00228.
This research was funded in whole or in part by National Science Centre (2021/43/I/ST3/00228). For the purpose of Open Access,
the author has applied a CC-BY public copyright licence to any Author Accepted Manuscript (AAM) version arising from this submission.

\appendix
\section{Proof of error estimate for $m < 2 \alpha$}
\label{sec:appendixErrorProof}
For convenience let us concisely repeat the proposition from the main text:
\newtheorem{prop}{Proposition}
\begin{prop}
Let $u \in W^{2 \alpha+1}_\infty(\Omega) \cap C^{2\alpha}(\Omega)$ and assume its polyharmonic spline interpolant $I_h u$ exists. Additionally assume that the order of the polyharmonic splines is $k > 2 \alpha$ and the approximation is augmented with monomials of degree $m < 2\alpha$. For sufficiently small $h$ it holds that 
    \begin{align}
    \norm{ \Delta^{\alpha}(I_h u - u)}_{L^\infty(K_i)} &\le C_i h^{m+1 -2\alpha} \abs{u}_{W^{m+1}_\infty(K_i)}.
    \label{eq:invers_inequality}
    \end{align}
\end{prop}
\newproof{proof}{Proof}
\begin{proof}
We observe that polyharmonic splines are in $C^{k-1}({K_i})$, hence $I_h u \in C^{2\alpha}({K_i})$, since $k > 2\alpha$. We now consider an auxiliary interpolation operator $I_h': C(\Omega) \to V_h'$, where $V_h'$ is obtained by augmentation of degree $m^\prime = 2 \alpha$. By the canonical splitting
    \begin{align}
        \Delta^\alpha I_h u - \Delta^\alpha u = \Delta^\alpha I_h u - \Delta^\alpha I_h' u  + \Delta^\alpha I_h' u - \Delta^\alpha u
        \label{eq:splitting}
    \end{align}
    Notice that $I_hu - I_h' u$ lies on a finite dimensional space $V_h + V_h' \subset W^{2\alpha}_\infty(K_i)$. This, combined with the well-known scale invariance arguments and the fact that our Voronoi regions are shape-regular, allows us to use the inverse inequalities~\cite[Lemma (4.5.3)]{Brenner2008} to estimate
    \begin{align}
        \|\Delta^\alpha I_h u - \Delta^\alpha I_h' u\|_{L^\infty(K_i)} \le C_1h^{-2\alpha} \| I_h u - I_h' u \|_{L^\infty(K_i)}. 
    \end{align}
    The factor appearing on the right-hand side can be further bound as
    \begin{align}
        \| I_hu - I'_h u \|_{L^\infty(K_i)} &\le \| I_hu - u \|_{L^\infty(K_i)} + \| I_h'u - u \|_{L^\infty(K_i)}\\
        & \le C_2h^{m+1} \abs{u}_{W^{m+1}_\infty(K_i)} + C_3h^{m'+1} \abs{u}_{W^{m'+1}_\infty(K_i)},
    \end{align}
    where the $h^{m+1}$ term dominates for small $h$, since $m^\prime > m$. Putting it all together we get
    \begin{align}
         \|\Delta^\alpha I_h u - \Delta^\alpha I_h' u\|_{L^\infty(K_i)} \le C_4h^{-2\alpha + m + 1} \abs{u}_{W^{m+1}_\infty(K_i)}.
    \end{align}
    To estimate the other term in equation~\eqref{eq:splitting}, we use the $m^\prime \ge 2 \alpha$ result to get
    \begin{align}
        \| \Delta^\alpha I_h' u - \Delta^\alpha u \|_{L^\infty(K_i)} \le  C_5 h^{m'+1 -2\alpha} \abs{u}_{W^{m^\prime+1}_\infty(K_i)}.
    \end{align}
    Again, the $h^{m+1}$ term dominates, which concludes the proof:
\begin{equation}
\begin{aligned}
    \norm{ \Delta^{\alpha}(I_h u - u)}_{L^\infty(K_i)} &\le C_6 h^{m+1 - 2\alpha} \abs{u}_{W^{m+1}_\infty(K_i)}.
\end{aligned}
\end{equation}
\end{proof}

\bibliographystyle{elsarticle-num}
\bibliography{adaptiveHyperviscosity}

\end{document}